\documentclass[11pt]{amsart}
\usepackage{indentfirst}
\usepackage{graphics,amsmath,amssymb,amsthm,mathrsfs,bm}
\usepackage{enumitem}

\usepackage{amsfonts}
\usepackage{amsmath}
\usepackage{amssymb}
\usepackage{multicol}
\usepackage{stmaryrd}
\usepackage{bbm}
\usepackage{cite}
\usepackage{epsfig}
\usepackage{color}
\usepackage{graphics}
\usepackage{graphicx}
\usepackage{epstopdf}
\usepackage{multicol,graphics}
\newcommand\bes{\begin{eqnarray}}
\newcommand\ees{\end{eqnarray}}
\newcommand\R{\mathbb R}
\newtheorem{theorem}{Theorem}[section]
\newtheorem{lemma}[theorem]{Lemma}
\newtheorem{corollary}[theorem]{Corollary}

\newtheorem{remark}[theorem]{Remark}

\numberwithin{equation}{section}
\allowdisplaybreaks

\begin{document}

\title[a degenerate epidemic model]{\textbf{The dynamics of a degenerate  epidemic model with nonlocal diffusion and free boundaries}}

\author[Zhao, Zhang, Li and Du]{Meng Zhao$^\dag$, Yang Zhang$^\ddag$, Wan-Tong Li$^{\dag}$ and Yihong Du$^{\S}$}

\thanks{\hspace{-.6cm}
$^\dag$ School of Mathematics and Statistics, Lanzhou University,
Lanzhou, Gansu 730000, P.R. China\\
$^\ddag$ Department of Mathematics, Harbin Engineering University, Harbin, 150001, P.R. China\\
$^\S$ School of Science and Technology, University of New England, Armidale, NSW 2351, Australia.\\
$^*${\sf Corresponding author} (wtli@lzue.edu.cn)
}

\date{\today}

\maketitle

\begin{abstract}
We consider an epidemic model with  nonlocal diffusion and free boundaries, which describes the evolution of  an infectious agents with nonlocal diffusion and the infected humans without diffusion, where humans get infected by the agents, and infected humans in return contribute to the growth  of the agents. The model can be viewed as a nonlocal version of the free boundary model  studied by Ahn, Beak and Lin \cite{ABL2016}, with its origin tracing back to Capasso et al. \cite{CP1979, CM1981}.
 We prove that the problem has a unique solution defined for all $t>0$, and its long-time dynamical behaviour is governed by a spreading-vanishing dichotomy. Sharp criteria for spreading and vanishing are also obtained, which reveal significant differences from the local diffusion model in \cite{ABL2016}. Depending on the choice of the kernel function in the nonlocal diffusion operator, it is expected that the nonlocal model here may have accelerated spreading, which would contrast sharply to the  model of \cite{ABL2016}, where the spreading has finite speed whenever spreading happens \cite{ZLN2019}.

\medskip


\smallskip


\end{abstract}


\section{Introduction}
\noindent

The spatial spread of epidemic disease is an important subject in mathematical epidemiology. In this paper we consider an epidemic model that
describes the evolution of  an infectious agents and the infected humans, where humans get infected by the agents, and infected humans in return contribute to the growth  of the agents.
In the model, the spatial movement of the infectious agents is described by a nonlocal diffusion operator, while that for the infective humans is ignored. The range of the infected area is assumed to be a moving interval  $[g (t), h(t)]\subset \mathbb R$, with its two end points representing the spreading fronts of the disease. Thus the model is a degenerate nonlocal diffusion system with free boundaries. We will show that the model has a unique solution defined for all time, and then determine  its long-time dynamical behaviour.

The origin of the model is the following ODE system
\begin{equation}\label{ODE}
\begin{cases}
u'(t)=-au(t)+cv(t),&t>0,\\
v'(t)=-bv(t)+G(u(t)),&t>0,
\end{cases}
\end{equation}
proposed by Capasso and Paveri-Fontana \cite{CP1979}
 to describe the cholera epidemic which spread in the European Mediterranean regions in 1973. Here   $a,\ b,\ c$ are all positive constants, $u(t)$ and $v(t)$, respectively, stand for the average population concentration of the infectious agents and the infective humans in the infected area at time $t$. $1/a$ stands for the mean lifetime of the agents in the environment,  $1/b$ represents the mean infectious period of the infective humans,  $c$ is the multiplicative factor of the infectious agents due to  the infective humans, and $G(u)$ is the infection rate of the human population due to the concentration of $u$ in the infected area. A basic assumption  of the model is that the  total susceptible human population is large enough compared  to the infective population, and is assumed to be constant during the evolution of the epidemic.

In \eqref{ODE}, the spatial factor is ignored. The corresponding spatial diffusion problem was subsequently considered by Capasso and Maddalena \cite{CM1981}, where it is assumed  that the infectious agents disperse randomly, and the mobility of the infective human population is small and thus neglected.  This diffusive model has the form
\begin{equation}\label{CM}
\begin{cases}
u_{t}=d\Delta u-au+cv,&t>0,\ x\in\Omega,\\
v_{t}=-bv+G(u),&t>0,\ x\in\Omega,\\
\frac{\partial u}{\partial \textbf{\emph{n}}}+\alpha u=0,
&t>0,\ x\in\partial\Omega,\\
u(0,x)=u_{0}(x),\ v(0,x)=v_{0}(x),&x\in\overline{\Omega},
\end{cases}
\end{equation}
where $\Omega\subset\mathbb R^N$ is a smooth bounded domain, representing the epidemic region, and the function $G$ is assumed to satisfy
\begin{enumerate}[leftmargin=2.8em]
\item[\textbf{(G1):}] $G\in C^{1}([0,\infty)),\ G(0)=0,\ G'(z)>0$ for $\forall z\geq 0$;
\item[\textbf{(G2):}] $\frac{G(z)}{z}$ is decreasing and $\lim\limits_{z\rightarrow +\infty}
    \frac{G(z)}{z}<\frac{ab}{c}$.
\end{enumerate}
It is shown in \cite{CM1981} that the number $$\tilde R_0:=\frac{cG'(0)}{(a+d\lambda_1)b}$$
is a threshold value for the long-time dynamical behaviour of \eqref{CM}:
The epidemic will eventually tend to extinction if $0<\tilde R_0\leq 1$, and there is a globally asymptotically stable endemic state  if $\tilde R_0>1$, where $\lambda_1$ is the first eigenvalue of
\[-\Delta\phi=\lambda\phi \text{\ in\ } \Omega,\;\;
\frac{\partial\phi}{\partial\textbf{\emph{n}}}+\alpha\phi=0
\text{\ on\ } \partial\Omega.\]

To describe how the epidemic spreads in space, one useful notion is the spreading speed. This can be achieved by considering \eqref{CM} over the entire $\mathbb R^N$ instead of over a bounded domain $\Omega$, coupled with initial functions $(u,v)=(u_0, v_0)$
which are positive over a bounded region, representing the infected area of the disease in the initial stage. A spreading speed can be established for this model, which is the minimal speed of its traveling wave solutions; we refer to Zhao and Wang \cite{ZW2004}, Wu at al. \cite{WSL2013} and references therein for research in this direction.

However, the approach described in the previous paragraph does not give the precise spreading front of the disease. This shortcoming can be addressed by considering the equations over a moving domain, resulting in a diffusive system with free boundaries.  Such an approach was taken by Ahn et al. \cite{ABL2016}, who considered the following free boundary version of  \eqref{CM} (in one space dimension),
\begin{equation}\label{FBlocal}
\begin{cases}
u_{t}=du_{xx}-au+cv,&t>0,\ x\in(g(t),h(t)),\\
v_{t}=-bv+G(u),&t>0,\ x\in(g(t),h(t)),\\
u(t,x)=v(t,x)=0,&t>0,\ x=g(t) \text{\ or\ } x=h(t),\\
g(0)=-h_{0},\ g'(t)=-\mu u_{x}(t,g(t)),&t>0,\\
h(0)=h_{0},\ h'(t)=-\mu u_{x}(t,h(t)),&t>0,\\
u(0,x)=u_{0}(x),\ v(0,x)=v_{0}(x),&x\in[-h_{0},h_{0}],
\end{cases}
\end{equation}
and proved a spreading-vanishing dichotomy for its long-time dynamical behaviour: The unique solution $(u,v,g,h)$ of \eqref{FBlocal}
satisfies one of the following:
\begin{itemize}
\item[{\rm(i)}] Vanishing:
\[
\mbox{$\lim\limits_{t\to\infty}[h(t)-g(t)]<\infty \text{\ and\ } \lim\limits_{t\rightarrow\infty} (\|u(t,\cdot)\|_{C([g(t),h(t)])}+\|v(t,\cdot)\|_{C([g(t),h(t)])})=0$;}
\]
\item[{\rm(ii)}] Spreading:
\[\mbox{$\lim\limits_{t\to\infty}[h(t)-g(t)]=\infty \text{\ and\ } R_0>1,\; \lim\limits_{t\rightarrow\infty}(u,v)=(K_1,K_2)
    \text{\ locally uniformly in\ } \mathbb{R}$,}
    \]
\end{itemize}
where
\begin{equation}\label{R0}
R_0:=\frac{cG'(0)}{ab},
\end{equation}
and $(K_1, K_2)$ are uniquely determined by
\[
G(K_1)=\frac{ab}{c}K_1, \;\; K_2=\frac{G(K_1)}{b}.
\]
Furthermore,  \begin{itemize}
\item[(i)] if $R_0\leq1$, then vanishing happens;
\item[(ii)] if $R_0\geq 1+ \frac da\big(\frac{\pi}{2h_0}\big)^2$, then spreading happens;
\item[(iii)]  if $1<R_0<1 +\frac da\big(\frac{\pi}{2h_0}\big)^2$, then vanishing happens  for  small initial data $(u_0, v_0)$, and spreading happens for  large initial data.
\end{itemize}
When spreading happens, the spreading speed of \eqref{FBlocal} was established in \cite{ZLN2019}.
\smallskip

Note that in \eqref{CM} and \eqref{FBlocal}, the dispersal of the infectious agents is assumed to follow the rules of random diffusion,
which is not realistic in general.  This kind of dispersal may be better described by a nonlocal diffusion operator of the form
\[d\displaystyle\int_{\R}J(x-y)u(t,y)dy-du(t,x),\]
which can capture short-range as well as long-range factors in  the dispersal by choosing the kernel function $J$ properly \cite{F2010,L2010}.

 The following nonlocal version of \eqref{CM},
\begin{equation}\label{uv-nonlocal}
\begin{cases}
u_t=d\displaystyle\int_{\R}J(x-y)u(t,y)dy-du-au+cv,& x\in\mathbb R, \; t>0,\\
v_t=-bv+G(u), & x\in\mathbb R,\; t>0
\end{cases}
\end{equation}
and its variations have been considered in several recent works; see, for example, Wang at al. \cite{WLS2018},  Zhang et al. \cite{ZLW2016}, and the references therein.

However, the nonlocal version of \eqref{FBlocal} has not been considered so far. After the work of
 Du and Lin \cite{DL2010} for a logistic type local diffusion model, free boundary approaches to local diffusion problems similar to \eqref{FBlocal} have been investigated by many researchers recently;
 see, e.g., \cite{DDL2019, DL2014,DL2015,DWZ2017, DWZ2018, GKLZ2015,
GLZ2015,GW2012,KM2018, KY2016, KLZ2013, LLS2016, L2019, LZ2017, WZ2015,W2014,WZ2014,WZ2017,WZ2018} and the references therein. Extensions of these free boundary problems to their nonlocal diffusion counterparts have been slow due partly to the fact that the Stefan condition in these local diffusion models
does not readily extend to the nonlocal problems.

Recently,  Cao et al. \cite{CDLL2018} proposed a nonlocal version of the logistic model of \cite{DL2010}, and successfully extended
many basic results of \cite{DL2010} to the nonlocal model.
In this paper, following the approach of \cite{CDLL2018}, we propose and examine a nonlocal version of \eqref{FBlocal}, which has the form
\begin{equation}\label{FB}
\begin{cases}
u_{t}=d\displaystyle\int_{g(t)}^{h(t)}J(x-y)u(t,y)dy-du
-au+cv,&t>0,\ x\in(g(t),h(t)),\\
v_{t}=-bv+G(u),&t>0,\ x\in(g(t),h(t)),\\
u(t,x)=v(t,x)=0,&t>0,\ x=g(t) \text{\ or\ } x=h(t),\\
g'(t)=-\mu\displaystyle\int_{g(t)}^{h(t)}\displaystyle\int_{-\infty}^{g(t)}
J(x-y)u(t,x)dydx,&t>0,\\
h'(t)=\mu\displaystyle\int_{g(t)}^{h(t)}\displaystyle\int_{h(t)}^{+\infty}
J(x-y)u(t,x)dydx,&t>0,\\
-g(0)=h(0)=h_{0},\ u(0,x)=u_{0}(x),\ v(0,x)=v_{0}(x),&x\in[-h_{0},h_{0}].
\end{cases}
\end{equation}
Here the kernel function $J: \R\rightarrow\R$ is assumed to satisfy
\begin{enumerate}[leftmargin=2em]
\item[\textbf{(J):}] $J\in C(\mathbb R)\cap L^\infty(\mathbb R)$,\; $J$ is symmetric and nonnegative,
\; $J(0)>0,\ \displaystyle\int_{\R}J(x)dx=1$. \end{enumerate}
 The parameters $a,\ b,\ c,\ d,\ \mu$ and $h_0$ are positive constants. The initial functions $u_0(x),v_0(x)$ satisfy
\begin{equation}\label{initial}
u_0,v_0\in C([-h_0,h_0]),\ u_0(\pm h_0)=v_0(\pm h_0)=0,
\ u_0,v_0>0 \text{\ in\ } (-h_0,h_0).
\end{equation}
As before, we assume $G$ satisfies {\bf (G1)-(G2)}. An example is $G(z)=\alpha\frac{z}{1+z}$ with $\alpha\in (0, ab/c)$.

In \eqref{FB}, the free boundary conditions
\[
\begin{cases}
\displaystyle h'(t)=\mu\int_{g(t)}^{h(t)}\int_{h(t)}^{+\infty}
J(x-y)u(t,x)dydx,\smallskip\\
\displaystyle g'(t)=-\mu\int_{g(t)}^{h(t)}\int_{-\infty}^{g(t)}
J(x-y)u(t,x)dydx,
\end{cases}
\]
mean that  the expanding rate of the range $[g(t), h(t)]$
is proportional to the outward flux of the population across the boundary of the range (see \cite{CDLL2018} for further explanations and justification).

\medskip

The main results of this paper are the following theorems:

\begin{theorem}[Global existence and uniqueness]\label{global-e-u}
Suppose that {\bf  (J)} and {\bf  (G1)-(G2)} hold. Then for any given $h_0>0$ and $u_0(x),v_0(x)$ satisfying \eqref{initial}, problem \eqref{FB} admits a unique solution $(u(t,x),v(t,x),g(t),h(t))$ defined for all $t>0$.
\end{theorem}

\begin{theorem}[Spreading-vanishing dichotomy]\label{dichotomy}
 Let the conditions of Theorem \ref{global-e-u} hold and $(u,v,g,h)$ be the unique solution of \eqref{FB}. Then one of the following alternatives must happen:
\begin{itemize}
\item[{\rm(i)}] {\rm Spreading:} $\lim\limits_{t\rightarrow\infty} [h(t)-g(t)]=\infty$  $($and  necessarily $R_0>1$$)$,
\[\lim\limits_{t\rightarrow+\infty}(u(t,x),v(t,x))=(K_1,K_2)
\text{\ locally uniformly in\ } \R.\]
\item[{\rm(ii)}] {\rm Vanishing:} $\lim\limits_{t\rightarrow\infty}
(g(t),h(t))=(g_\infty, h_\infty)$ is a finite interval,  and
\[\lim\limits_{t\rightarrow\infty}
\max\limits_{g(t)\leq x\leq h(t)}u(t,x)=0 \text{\ and\ } \lim\limits_{t\rightarrow\infty}
\max\limits_{g(t)\leq x\leq h(t)}v(t,x)=0.\]
\end{itemize}
\end{theorem}

Let us recall that $R_0$ is given by \eqref{R0}.

\begin{theorem}[Spreading-vanishing criteria]\label{criteria} In Theorem \ref{dichotomy}, the dichotomy can be determined as follows:
\begin{itemize}
\item[{\rm(i)}] If $R_0\leq 1$, then vanishing happens.
\item[{\rm(ii)}] If $R_0\geq 1+\frac da$, then spreading happens.
\item[{\rm(ii)}] If $1<R_0<1+\frac da$, then there exists $l^*>0$ such that spreading happens when $2h_{0}\geq l^\ast$, and if $2h_{0}<l^\ast$, then there exists  $\mu^*>0$ such that spreading happens if and only if  $\mu>\mu^*$.
\end{itemize}
\end{theorem}

We note that $l^*$ depends only on $(a,b,c, d, J)$, which is determined by an eigenvalue problem (see \eqref{eigen}), but $\mu^*$ depends also on the initial data.

\begin{remark}\label{rm1}{\rm
From part (ii) of Theorem 1.3, we see that if $R_0>1$ then for all small $d>0$, spreading happens. This is very different from the local diffusion model \eqref{FBlocal}, where in the corresponding case, the size of the initial population range $2h_0$ (and the initial functions) also plays an important role.}
\end{remark}

\begin{remark}\label{rm2}{\rm
Very recently, Du, Li and Zhou \cite{DLZ2019} investigated the spreading speed of the nonlocal model in \cite{CDLL2018} and proved that the spreading may or may not have a finite speed, depending on whether a certain condition is satisfied by the kernel function $J$ in the nonlocal diffusion term. We expect a similar result for \eqref{FB}, which will be considered in a future work.}
\end{remark}

In \cite{DWZ2019}, some two species Lotka-Volterra models with nonlocal diffusion and free boundaries have been considered. There, nonlocal diffusion happens to both species, and the reaction/growth functions are also very different from \eqref{FB} here. As a result, the techniques and results there are very different from this paper here.

The rest of the paper is organised as follows. In Section 2 we prove Theorem \ref{global-e-u}, namely, problem \eqref{FB} has a unique solution defined for all $t>0$. The long-time dynamical behaviour of \eqref{FB} is investigated in Section 3,  where Theorems \ref{dichotomy} and \ref{criteria}  are proved.

\section{Global existence and uniqueness}
\noindent

Throughout this section, we assume that $h_0>0$  and $(u_0, v_0)$ satisfy \eqref{initial}. For any given $T>0$, we  introduce the following notations:
\begin{align*}
&A:=\max\left\{K_1,\ \|u_0\|_\infty,\ \frac{c}{a}\|v_0\|_\infty\right\},\\
&B:=\max\left\{\|v_0\|_\infty,\ \frac{G(A)}{b}\right\},\\
&H_T=H^{h_0}_{T}:=\left\{h\in C([0,T])\ :\ h(0)=h_0,\
\inf_{0\leq t_1<t_2\leq T}
\frac{h(t_2)-h(t_1)}{t_2-t_1}>0\right\},\\
&G_T=G^{h_0}_{T}:=\{g\in C([0,T])\ :\ -g\in H^{h_0}_{T}\},\\
&D_T^{g,h}:=\left\{(t,x)\in\R^2\ :\ 0<t\leq T,\ g(t)<x<h(t)\right\},\\
&D_T=D_T^{h_0}:=\left\{(t,x)\in\R^2\ :\ 0<t\leq T,\ -h_0<x<h_0\right\},\\
&\mathbb{X}^{v_0}_T:=\left\{\phi\in C(\overline D_T^{g,h})\ :
\ \phi(0,x)=v_0(x) \text{\ in\ } [-h_0,h_0],
\ 0\leq\phi\leq B \text{\ in\ } D_T^{g,h},\right.\\
&\hspace{5cm} \phi(t,x)=0 \text{\ for\ } t\in(0,T],\ x\in\R\backslash(g(t),h(t))\Big\},\\
&\mathbb{X}^{u_0}_T:=\left\{\phi(t,x)\in C(\overline D_T^{g,h})\ :
\ \phi(0,x)=u_0(x) \text{\ in\ } [-h_0,h_0],\right.\\
&\hspace{5cm} \phi(t,x)=0 \text{\ for\ } t\in(0,T],\ x\in\R\backslash(g(t),h(t))\Big\}.
\end{align*}
Noting that $A\geq K_1$ and $\frac{G(K_1)}{K_1}=\frac{ab}{c}$, we have $G(A)\leq\frac{ab}{c}A$ by (G2), and hence
\begin{equation}\label{AB}
B\leq\frac{a}{c}A.
\end{equation}

The following maximum principle will be frequently used in our discussions below.

\begin{lemma}\label{MP}
Assume that {\bf  (J)} holds, and $g\in G_{T},\ h\in H_{T}$ for some $T>0$. Suppose that $c_{ij}\in L^\infty(D_T^{g,h})$
for $i, j\in \{1,2\}$ with $c_{12}$ and $c_{21}$ nonnegative, and $(u(t,x),v(t,x))$ as well as $(u_t(t,x),v_t(t,x))$ are continuous in $\overline{D_T^{g,h}}$ and satisfy
\begin{equation}\label{uv>=00}
\begin{cases}
u_{t}\geq d\displaystyle\int_{g(t)}^{h(t)}J(x-y)u(t,y)dy-du
+c_{11}u+c_{12}v,&0<t\leq T,\ x\in(g(t),h(t)),\\
v_{t}\geq c_{21}u+c_{22} v,&0<t\leq T,\ x\in(g(t),h(t)),\\
u(t,x)\geq0,\ v(t,x)\geq0,&0<t\leq T,\ x=g(t) \text{\ or\ } h(t),\\
u(0,x)\geq0,\ v(0,x)\geq0,&x\in[-h_{0},h_{0}].
\end{cases}
\end{equation}
Then $(u(t,x),v(t,x))\geq(0,0)$ for all $0\leq t\leq T$ and
$g(t)\leq x\leq h(t)$. Moreover, if we assume additionally $u(0,x)\not\equiv0$ in $[-h_0,h_0]$, then $u(t,x)>0$ in $D_T^{g,h}$.
\end{lemma}

\begin{proof}
Let $w(t,x)=e^{kt}u(t,x)$ and $z(t,x)=e^{kt}v(t,x)$, where $k$ is large enough such that
\[k>d+\|c_{11}\|_\infty+\|c_{22}\|_\infty,\]
and then
\[p(t,x):=k+c_{22}(t,x)+c_{21}(t,x)\geq k-\|c_{22}\|_\infty>0 \text{\ for all\ } (t,x)\in D_T^{g,h}.\]
By direct calculations, we have that $w(t,x)$ and $z(t,x)$ satisfy
\begin{equation}\label{wz--}
\begin{cases}
w_{t}\geq d\displaystyle\int_{g(t)}^{h(t)}J(x-y)w(t,y)dy+(k-d+c_{11})w+c_{12}z,&\\
z_{t}\geq c_{21}w+(k+c_{22})z.&
\end{cases}
\end{equation}
Denote
\[
p_0=\sup\limits_{(t,x)\in D_T^{g,h}}p(t,x) \; \mbox{  and } \;T^\ast=\min\left\{T,\ \frac{1}{4(k+\|c_{11}\|_\infty+\|c_{12}\|_\infty)},\ \frac{1}{4p_0}\right\}.
\]
 Now we prove $w\geq0$ and $z\geq0$ in $D_{T^\ast}^{g,h}$. Suppose that \[m:=\min\left\{\inf\limits_{(t,x)\in D_{T^\ast}^{g,h}}w(t,x),
\inf\limits_{(t,x)\in D_{T^\ast}^{g,h}}z(t,x)\right\}<0.\]
By \eqref{uv>=00}, $w\geq0$ and $z\geq0$ on the boundary of $D_{T^\ast}^{g,h}$. Hence, there exists $(t^\ast,x^\ast)\in D_{T^\ast}^{g,h}$ such that $\frac{m}{2}=w(t^\ast,x^\ast)<0$ or $\frac{m}{2}=z(t^\ast,x^\ast)<0$. We now define
\begin{equation*}
t_0=t_0(x^\ast):=
\begin{cases}
t_{x^\ast}^{g},&x^\ast\in(g(t^\ast),-h_{0}) \text{\ and\ } x^\ast=g(t_{x^\ast}^{g}),\\
0,&x^\ast\in[-h_{0},h_{0}],\\
t_{x^\ast}^{h},&x^\ast\in(h_{0},h(t^\ast)) \text{\ and\ }
x^\ast=h(t_{x^\ast}^{h}).
\end{cases}
\end{equation*}
Clearly, $u(t_0,x^\ast)\geq0$ and $v(t_0,x^\ast)\geq0$.

If $\frac{m}{2}=w(t^\ast,x^\ast)<0$, then it follows from the choice of $k$ and the first equation of \eqref{wz--} that
\begin{align*}
w(t^\ast,x^\ast)-w(t_0,x^\ast)\geq&\ d\displaystyle\int_{t_0}^{t^\ast}
\displaystyle\int_{g(t)}^{h(t)}J(x^\ast-y)w(t,y)dydt\\
\quad &+\displaystyle\int_{t_0}^{t^\ast}[(k-d+c_{11})w(t,x^\ast)+c_{12}z(t,x^\ast)]dt\\
\geq&\ d\displaystyle\int_{t_0}^{t^\ast}
\displaystyle\int_{g(t)}^{h(t)}J(x^\ast-y)mdydt
+\displaystyle\int_{t_0}^{t^\ast}[(k-d+c_{11})m+c_{12}m]dt\\
\geq&\ m(k+\|c_{11}\|_\infty+\|c_{12}\|_\infty )(t^\ast-t_0).
\end{align*}
Since $w(t_0,x^\ast)=e^{kt_0}u(t_0,x^\ast)\geq0$, we deduce
\[\frac{m}{2}\geq m(k+\|c_{11}\|_\infty+\|c_{12}\|_\infty)(t^\ast-t_0)\geq m(k+\|c_{11}\|_\infty+\|c_{12}\|_\infty)T^\ast\geq\frac{m}{4},\]
which is a contradiction to $m<0$.

If $\frac{m}{2}=z(t^\ast,x^\ast)<0$, then it follows from the choice of $k$ and the second equation of \eqref{wz--} that
\begin{align*}
z(t^\ast,x^\ast)-z(t_0,x^\ast)\geq&\
\displaystyle\int_{t_0}^{t^\ast}[(k+c_{22})z(t,x^\ast)
+c_{21}w(t,x^\ast)]dt\\
\geq&\ \displaystyle\int_{t_0}^{t^\ast}[(k+c_{22})m
+c_{21}m]dt\\
\geq&\ mp_0(t^\ast-t_0).
\end{align*}
Since $z(t_0,x^\ast)=e^{kt_0}v(t_0,x^\ast)\geq0$, we deduce
\[\frac{m}{2}\geq mp_0(t^\ast-t_0)\geq mp_0T^\ast\geq\frac{m}{4},\]
which is a contradiction to $m<0$.

If $T^\ast=T$, then $(u(t,x),v(t,x))\geq(0,0)$ for all $0\leq t\leq T$ and $g(t)\leq x\leq h(t)$ follows directly; while if $T^\ast<T$, we may repeat this process with $(u_0(x),v_0(x))$ replaced by $(u(T^\ast,x),v(T^\ast,x))$, and $(0,T]$ replaced by $(T^\ast,T]$. Clearly after repeating this process finitely many times, we will obtain $(u(t,x),v(t,x))\geq(0,0)$ for all $0\leq t\leq T$ and $g(t)\leq x\leq h(t)$.

Due to $v(t,x)\geq0$ for all $0\leq t\leq T$ and $g(t)\leq x\leq h(t)$, we have that $u$ satisfies
\begin{equation*}
\begin{cases}
u_t\geq d\displaystyle\int_{g(t)}^{h(t)}J(x-y)u(t,y)dy-du
+c_{11}u,&0<t\leq T,\ x\in(g(t),h(t)),\\
u(t,g(t))\geq0,\ u(t,h(t))\geq0,&0<t\leq T,\\
u(0,x)\geq0,&x\in[-h_0,h_0].
\end{cases}
\end{equation*}
If $u(0,x)\not\equiv0$ in $[-h_0,h_0]$, then it follows directly from Lemma 2.2 of \cite{CDLL2018} that $u(t,x)>0$ in $D_T^{g,h}$.
\end{proof}

The following result is an important first step towards the proof of Theorem \ref{global-e-u}.

\begin{lemma}\label{ugh-e-u}
For any given $T>0$, $(g,h)\in G_{T}\times H_{T}$ and $v\in\mathbb{X}^{v_0}_T$, the problem
\begin{equation}\label{uODE}
\begin{cases}
u_t=d\displaystyle\int_{g(t)}^{h(t)}J(x-y)u(t,y)dy-du
-au+cv,&0<t\leq T,\ x\in(g(t),h(t)),\\
u(t,g(t))=u(t,h(t))=0,&0<t\leq T,\\
u(0,x)=u_0(x),&x\in[-h_0,h_0]
\end{cases}
\end{equation}
admits a unique solution  $u^\ast\in C(\overline{D_T^{g,h}})$.
Moreover,
\begin{equation}\label{uA}
0<u^\ast\leq A \text{\ for any\ } (t,x)\in D_T^{g,h}.
\end{equation}
\end{lemma}
\begin{proof}
For any given $v\in \mathbb{X}^{v_0}_T$, let $f(t,x,u)=-au+cv(t,x)$. Since $f(t,x,0)\not\equiv 0$, the corresponding result in \cite{CDLL2018}
does not cover the case here. However, the method in \cite{CDLL2018} can be extended to deal with this case. Since considerable changes are needed, we give the details below for completeness.

\textbf{Step 1}: A parameterised ODE problem.

For any given $x\in[g(T),h(T)]$, define
\begin{equation*}
\widehat{u}_{0}(x):=
\begin{cases}
u_{0}(x),&x\in[-h_{0},h_{0}],\\
0,&x\not\in[-h_{0},h_{0}]
\end{cases}
\end{equation*}
and
\begin{equation*}
t_{x}:=
\begin{cases}
t_{x}^{g},&x\in[g(T),-h_{0}) \text{\ and\ } x=g(t_{x}^{g}),\\
0,&x\in[-h_{0},h_{0}],\\
t_{x}^{h},&x\in(h_{0},h(T)] \text{\ and\ }
x=h(t_{x}^{h}).
\end{cases}
\end{equation*}
Clearly, $t_x=T$ for $x=g(T)$ or $x=h(T)$, and $0\leq t_x<T$ for $x\in(g(T),h(T))$.

For any given $s\in(0,T]$ and $\phi\in\mathbb{X}^{u_0}_s$, we fix $x\in (g(s), h(s))$ and consider the  problem
\begin{equation}\label{uphiODE}
\begin{cases}
u_{t}(t,x)=F(t,x,u),&t_{x}<t\leq s,\\
u(t_{x},x)=\widehat{u}_{0}(x),&
\end{cases}
\end{equation}
with
\[
F(t,x,u):=d\displaystyle\int_{g(t)}^{h(t)}J(x-y)\phi(t,y)dy-du-au+cv(t,x).
\]
We will regard \eqref{uphiODE} as an ODE initial value problem with parameter $x$.
Set
 \[L_{1}:=1+\max\left\{\|\phi\|_{C\left(\overline{D_T^{g,h}}\right)},
A\right\}.\]
For any $u_1,u_2\in[0,L_{1}]$,
\[|F(t,x,u_{1})-F(t,x,u_{2})|=(d+a)|u_{1}-u_{2}|,\]
Hence, $F(t,x,u)$ is Lipschitz continuous in $u$ for $u\in[0,L_1]$ with Lipschitz constant $d+a$, uniformly for $t\in[0,s]$ and $x\in(g(s),h(s))$. Additionally, $F(t,x,u)$ is continuous in all its variables in this range. By the fundamental theorem of ODEs, problem \eqref{uphiODE} admits a unique solution $U_\phi(t,x)$ defined in some interval $[t_{x},s_x)$ of $t$, and $U_\phi(t,x)$ is continuous in both $t$ and $x$.

To see that $t\rightarrow U_\phi(\cdot, x)$ can be uniquely extended to $[t_{x},s]$, it suffices to show that if $U_\phi$ is uniquely defined for $t\in[t_{x},\widehat{t}]$ with $\widehat{t}\in(t_{x},s]$, then
\begin{equation}\label{U-phi-L1}
0\leq U_\phi(t,x)\leq L_{1} \text{\ for\ } t\in[t_{x},\widehat{t}].
\end{equation}
 It is easy to check that
\begin{align*}
F(t,x, L_1)=d\displaystyle\int_{g(t)}^{h(t)}J(x-y)\phi(t,y)dy-dL_1-aL_1
+cv\leq d\|\phi\|_\infty-dL_1-aA+cB\leq0
\end{align*}
and
\[L_1>\|\phi\|_{C\left(\overline{D_T^{g,h}}\right)}
\geq\|u_0(x)\|_{C([-h_0,h_0])}=\|\widehat{u}_0\|_\infty.\]
Now a simple comparison argument gives  $U_\phi(t,x)\leq L_{1}$ for $t\in[t_{x},\widehat{t}]$. This proves the second inequality in \eqref{U-phi-L1}. The first inequality there can be obtained similarly by using $F(t,x,0)\geq 0$.

\textbf{Step 2}: A fixed point problem.

For $s\in (0, T]$, for simplicity we denote
\[D_s:=D_s^{g,h},\ \mathbb{X}_s:=\mathbb{X}^{u_0}_s.\]
 By Step 1, for any $\phi\in\mathbb{X}_{s}$ we can find a unique $U_\phi(t,x)$ satisfying \eqref{uphiODE} for $t\in [0, s]$, and  by the continuous dependence of the ODE solution on parameters, $U_\phi(t,x)$ is continuous in $\overline D_s$. Hence $U_\phi\in\mathbb{X}_{s}$. Note that $\mathbb{X}_s$ is a complete metric space equipped with the norm
\[d(\phi_1,\phi_2)=\|\phi_1-\phi_2\|_{C(\overline{D}_s)}.\]
Now, we define a mapping $\Gamma:\mathbb{X}_s\longrightarrow
\mathbb{X}_s$ by
\[\Gamma(\phi)=U_\phi.\]
Clearly, if $\Gamma(\phi)=\phi$, then $\phi$ solves \eqref{uODE} for $t\in [0, s]$, and vice versa.

Letting $M=\max\left\{2\|u_{0}\|_\infty,A\right\}$, we denote
\[\mathbb{X}_s^M:=\left\{\phi\ |\ \phi\in\mathbb{X}_s,\ \phi(t,x)\leq M \mbox{ in } D_s\right\}.\]
In the following, we will show that $\Gamma$ has a unique fixed point in $\mathbb{X}_s^M$ for small $s$ by the contraction mapping theorem.

We first claim that there exists sufficiently small $s^{\ast}$ such that $\Gamma$ maps $\mathbb{X}^M_s$ into itself for any $s\in(0,s^{\ast}]$. So let $\phi\in\mathbb{X}^M_s$, and we aim to show that $U_\phi(t,x)\leq M$ for $(t,x)\in\overline{D}_s$. By the first equation of \eqref{uphiODE} and \eqref{AB}, we have
\[(U_\phi)_{t}(t,x)\leq d\displaystyle\int_{g(t)}^{h(t)}J(x-y)\phi(t,y)dy+cB\leq d\|\phi\|_{C(\overline{D}_s)}+aA \text{\ for\ } t_{x}<t\leq s,\]
and then
\begin{align*}
U_\phi(t,x)
\leq\ &U_\phi(t_{x},x)
+(d\|\phi\|_{C(\overline{D}_s)}+aA)(t-t_x)\\
\leq\ &\|u_{0}\|_{\infty}
+(d\|\phi\|_{C(\overline{D}_s)}+aA)s\\
\leq\ &\frac{M}{2}+(d+a)sM.
\end{align*}
If we choose $s^{\ast}$ small enough such that
\[(d+a)s^\ast\leq\frac{1}{2},\]
then $U_\phi(t,x)\leq M$ for $s\in(0,s^\ast]$. This implies that $U_\phi\in\mathbb{X}_s^M$ for $s\in(0,s^\ast]$. The claim is now proved.

Next, we show that $\Gamma$ is a contraction map for $s\in\left(0,s^\ast\right]$.  Namely,  there exists some $\delta<1$ such that for any given $\phi_{i}\in \mathbb X_s^M$, $ i=1,2$, we have
\[\|U_{\phi_1}-U_{\phi_2}\|_{C(\overline{D}_s)}
\leq\delta\|\phi_{1}-\phi_{2}\|_{C(\overline{D}_s)}.\]
For such $\phi_1$ and $\phi_2$, denote $W=U_{\phi_1}-U_{\phi_2}$. Then
\begin{equation*}
\begin{cases}
W_{t}+(d+a)W=d\displaystyle\int_{g(t)}^{h(t)}J(x-y)
(\phi_{1}-\phi_{2})(t,y)dy,&t_{x}<t\leq s,\\
W(t_{x},x)=0,&x\in(g(s),h(s)).
\end{cases}
\end{equation*}
By direct calculations, we have
\[W(t,x)=de^{-(d+a)t}\displaystyle\int_{t_{x}}^{t}
e^{(d+a)\tau}\displaystyle\int_{g(\tau)}^{h(\tau)}J(x-y)
(\phi_{1}-\phi_{2})(\tau,y)dyd\tau \text{\ for\ } t_{x}<t\leq s.\]
Therefore
\[|W(t,x)|\leq d\|\phi_{1}-\phi_{2}\|_{C(\overline{D}_s)}(t-t_x)
\leq ds\|\phi_{1}-\phi_{2}\|_{C(\overline{D}_s)}.\]
It follows that
\[\|W\|_{C(\overline{D}_s)}\leq\frac{1}{2}
\|\phi_{1}-\phi_{2}\|_{C(\overline{D}_s)} \text{\ for\ } s\in(0,s^\ast].\]
Hence, $\Gamma$ is a contraction map. For any $s\in(0,s^\ast]$, applying the contraction mapping theorem we obtain  a unique fixed point $u^\ast$ of $\Gamma$ in $\mathbb{X}_s^M$. Clearly $u^\ast$ solves \eqref{uODE} for $t\in [0, s]$.

To see $u^*$ is the unique solution of \eqref{uODE} for $t\in[0,s]$ with $s\in(0,s^\ast]$, it remains to show that any solution $u$ of \eqref{uODE} for $t\in [0, s]$ belongs to $\mathbb X_s^M$. We show below that actually the following sharper estimates hold:
\begin{equation}\label{UAM}
0\leq u(t,x)\leq A \text{\ for\ } t\in[0,s] \text{\ and\ } x\in[g(t),h(t)].
\end{equation}
It is easy to check that
\begin{align*}
&u_t\leq d\displaystyle\int_{g(t)}^{h(t)}J(x-y)u(t,y)dy-du
-au+cB,\\
& d\displaystyle\int_{g(t)}^{h(t)}J(x-y)Ady-dA-aA+cB\leq -aA+cB\leq 0.
\end{align*}
Therefore, in view of $A\geq\|u_0\|_\infty$, a simple comparison argument yields  $u(t,x)\leq A$ for $t\in[0,s]$ and $x\in[g(t),h(t)]$. Using
\begin{equation*}
\begin{cases}
u_t\geq d\displaystyle\int_{g(t)}^{h(t)}J(x-y)u(t,y)dy-du
-au,&0<t\leq s,\ x\in(g(t),h(t)),\\
u(t,g(t))=u(t,h(t))=0,&0<t\leq s,\\
u(0,x)=u_0(x)\geq0,&x\in[-h_0,h_0],
\end{cases}
\end{equation*}
we see from Lemma 2.2 of \cite{CDLL2018} that $u(t,x)\geq0$ for $t\in[0,s]$ and $x\in[g(t),h(t)]$. Hence, \eqref{UAM} holds.
We have thus proved that for any $s\in (0, s^*]$, \eqref{uODE} has a unique solution for $t\in [0, s]$.

\textbf{Step 3}: Extension of the solution.

 Since $s^*$ in Step 2 only depends on $a$ and $b$, we may repeat Step 2 to \eqref{uODE} with the initial time $t=0$ replaced by $t=s$
 for any $s\in (0, s^*]$, and so the solution of \eqref{uODE} can be uniquely extended to $t\in [0, s]$ for any $s\in (0, \min\{2s^*, T\}]$.
  Moreover, the extended solution $u$ still satisfies \eqref{UAM}.  By repeating this process finitely many times, the solution of problem \eqref{uODE} is uniquely extended to $[t_{x},T]$, and \eqref{uA} is a consequence of \eqref{UAM}
  obtained in each step of the extension.
\end{proof}

\noindent
{\bf Proof of Theorem \ref{global-e-u}:} Following the approach of \cite{CDLL2018}, we make use of Lemma \ref{ugh-e-u} and a fixed point argument.
 For any given $T>0$ and $(g^\ast,h^\ast)\in G_T\times H_T$, $v^\ast\in\mathbb{X}_T^{v_0}$, it follows from Lemma \ref{ugh-e-u} that \eqref{uODE} with $(v,g,h)=(v^\ast,g^\ast,h^\ast)$ has a unique solution $u^\ast$. For such $(u^\ast,g^\ast,h^\ast)$, we can define
 $\widehat v_0(x)$ as the zero extension of $v_0(x)$ to $x\in\mathbb R\setminus [-h_0,h_0]$ and then define $t_x$ as in Step 1 of the proof of Lemma \ref{ugh-e-u}, but with $(g,h)$ replaced by $(g^*, h^*)$. To mark the difference, we denote $t_x$ by $t^*_x$.

 Now, for each $x\in (g^*(T), h^*(T))$, we consider the initial value problem
\begin{equation}\label{vgh}
\begin{cases}
v_t=-bv+G(u^\ast),& t_x^*<t\leq T,\\
v(t_x^*,x)=\widehat v_0(x).&
\end{cases}
\end{equation}
By the Fundamental Theorem of ODEs and some simple comparison argument, it can be easily shown that \eqref{vgh} has a unique solution $\widetilde v^*(t,x)$, and it is continuous and satisfies
\[\mbox{
$0\leq \widetilde{v}^\ast(t,x)\leq B$ for $t\in(0,T]$ and $x\in[g^\ast(t),h^\ast(t)]$}.
\]
Therefore $\widetilde{v}^\ast\in\mathbb{X}^{v_0}_T$.

Next we define $(\widetilde{g}^\ast,\widetilde{h}^\ast)$ for $t\in [0, T]$ by
\begin{equation}\label{gh-wide}
\begin{cases}
&\widetilde{g}^\ast(t):=-h_0-\mu\displaystyle\int_{0}^{t}
\displaystyle\int_{g^\ast(\tau)}^{h^\ast(\tau)}
\displaystyle\int_{-\infty}^{g^\ast(\tau)}
J(x-y)u^\ast(\tau,x)dydxd\tau,\\
&\widetilde{h}^\ast(t):=h_0+\mu\displaystyle\int_{0}^{t}
\displaystyle\int_{g^\ast(\tau)}^{h^\ast(\tau)}
\displaystyle\int_{h^\ast(\tau)}^{+\infty}
J(x-y)u^\ast(\tau,x)dydxd\tau.
\end{cases}
\end{equation}
Due to {\bf (J)}, there exist constants $\epsilon_0\in(0,h_0/4)$ and $\delta_0$ such that
\[J(x)\geq\delta_0 \text{\ if\ } |x|\leq\epsilon_0.\]
Using this we can follow the corresponding arguments of \cite{CDLL2018} to show that, for
 some sufficiently small $T_0=T_0(\mu,A,h_0,\epsilon_0,u_0,J)>0$ and any $T\in (0, T_0]$,
\[
\sup\limits_{0\leq t_1<t_2\leq T}
\frac{\widetilde g^*(t_2)-\widetilde g^*(t_1)}{t_2-t_1}\leq
-\mu\widetilde{\sigma}_0,\ \ \
\inf\limits_{0\leq t_1<t_2\leq T}
\frac{\widetilde h^*(t_2)-\widetilde h^*(t_1)}{t_2-t_1}\geq
\mu\sigma_0,\
\]
\[
\widetilde h^*(t)-\widetilde g^*(t)\leq2h_0+\frac{\epsilon_0}{4}
\text{\ for\ } t\in[0,T],
\]
where
\[\widetilde{\sigma}_0=\frac{1}{4}\epsilon_0\delta_0e^{-(d+a)T_0}
\displaystyle\int_{-h_0}^{-h_0+\frac{\epsilon_0}{4}}u_0(x)dx,\;\;
\sigma_0=\frac{1}{4}\epsilon_0\delta_0e^{-(d+a)T_0}
\displaystyle\int_{h_0-\frac{\epsilon_0}{4}}^{h_0}u_0(x)dx.\]

Let
\begin{align*}
&\Sigma_T:=\left\{(v,g,h)\in \mathbb{X}^{v_0}_T\times G^{h_0}_{T}\times H^{h_0}_{T}\ :\ \sup\limits_{0\leq t_1<t_2\leq T}
\frac{g(t_2)-g(t_1)}{t_2-t_1}\leq
-\mu\widetilde{\sigma}_0,\right.\\
&\ \ \ \ \ \ \ \ \left.
\ \ \ \ \inf\limits_{0\leq t_1<t_2\leq T}
\frac{h(t_2)-h(t_1)}{t_2-t_1}\geq
\mu\sigma_0,\
h(t)-g(t)\leq2h_0+\frac{\epsilon_0}{4}
\text{\ for\ } t\in[0,T]\right\},
\end{align*}
and define the mapping
\[\mathcal{F}(v^\ast,g^\ast,h^\ast)
=(\widetilde{v}^\ast,\widetilde{g}^\ast,\widetilde{h}^\ast).\]
Then the above analysis indicates that
\[\mathcal{F}(\Sigma_T)\subset\Sigma_T \text{\ for\ } T\in(0,T_0].\]

In the following, we show that for sufficiently small $T\in (0, T_0]$, $\mathcal {F}$ has a unique fixed point in $\Sigma_T$, which clearly
is a solution of \eqref{FB} for $t\in [0, T]$. We will then show that this is the unique solution of \eqref{FB} and it can be extended uniquely to all $t>0$. We will complete this task in several steps.

\textbf{Step 1:} We show that, for sufficiently small $T\in(0,T_0]$, $\mathcal{F}$ has a unique fixed point in $\Sigma_T$ by the contraction mapping theorem.

For $T\in (0, T_0]$ and any given $(v^\ast_i,g^\ast_i,h^\ast_i)\in\Sigma_T\ (i=1,2)$, denote
\[
(\widetilde v_i^*, \widetilde{g}^\ast_i,\widetilde{h}^\ast_i)=\mathcal F (v^\ast_i,g^\ast_i,h^\ast_i),\;\; i=1,2.
\]
Define
\begin{align*}
&h_{m}(t):=\min\{h_{1}^\ast(t),\ h_{2}^\ast(t)\},
\ h_{M}(t):=\max\{h_{1}^\ast(t),\ h_{2}^\ast(t)\},\\
&g_{m}(t):=\min\{g_{1}^\ast(t),\ g_{2}^\ast(t)\},
\ g_{M}(t):=\max\{g_{1}^\ast(t),\ g_{2}^\ast(t)\}.
\end{align*}
In the following, we will show that there exists $\delta\in(0,1)$ such that, for all small $T\in (0, T_0]$, and any $(v^\ast_i,g^\ast_i,h^\ast_i)\in\Sigma_T\ (i=1,2)$,
\begin{equation}\label{v-wast-dt}
\begin{split}
&\|\widetilde{v}^\ast_1-\widetilde{v}^\ast_2\|_{C([0,T]\times\mathbb R )}
+\|\widetilde{g}^\ast_1-\widetilde{g}^\ast_2\|_{C([0,T])}
+\|\widetilde{h}^\ast_1-\widetilde{h}^\ast_2\|_{C([0,T])}\\
\leq\ &\delta\left(\|v^\ast_1-v^\ast_2\|_{C([0,T]\times\mathbb R)}
+\|g^\ast_1-g^\ast_2\|_{C([0,T])}
+\|h^\ast_1-h^\ast_2\|_{C([0,T])}\right).
\end{split}
\end{equation}
Clearly this implies that $\mathcal F$ is a contraction mapping on $\Sigma_T$.

To prove \eqref{v-wast-dt}, we first estimate $\|\widetilde{v}^\ast_1-\widetilde{v}^\ast_2\|_{C([0,T]\times\mathbb R)}$. Let
\begin{align*}
&\widetilde{V}(t,x):=
\widetilde{v}^\ast_1(t,x)-\widetilde{v}^\ast_2(t,x),\ \
 U(t,x):=u^\ast_1(t,x)-u^\ast_2(t,x),\ \
 V(t,x):=v^\ast_1(t,x)-v^\ast_2(t,x).
\end{align*}

\smallskip

{\bf Claim 1:} There exist positive constants $\tilde T$ and $C$ such that for any $T\in (0, \tilde T]$ and any $(t^\ast,x^\ast)\in D_T^{g_m,h_M}$,
\begin{equation}\label{tilde-V-*}
|\widetilde{V}(t^\ast,x^\ast)|\leq CT\left[\|V\|_{C([0,T]\times\R)}
+\|g_1^\ast-g_2^\ast\|_{C([0,T])}
+\|h_1^\ast-h_2^\ast\|_{C([0,T])}\right].
\end{equation}

To prove Claim 1,  we proceed according to three cases.
\smallskip

\underline{Case 1}: $x^\ast\in[-h_0,h_0]$.

From \eqref{vgh} we obtain
\begin{align*}
\widetilde V(t^*,x^*)
=\ &e^{-bt^*}\left[\displaystyle\int_0^{t^*}
e^{b\tau}\left(G(u^\ast_1)-G(u^\ast_2)\right)(\tau,x^*)d\tau\right].
\end{align*}
Since $G\in C^1([0,\infty))$, for any $L>0$, there exists a constant $\rho(L)>0$ such that
\[|G(z_1)-G(z_2)|\leq\rho(L)|z_1-z_2| \text{\ for\ } z_1,z_2\in[0,L].\]
It follows that
\begin{equation}\label{VI}
|\widetilde V(t^*,x^*)|
\leq e^{-bt^*}\displaystyle\int_0^{t^*}e^{bt^*}d\tau\rho(A)
\|U\|_{C([0,T]\times\R)}
\leq\rho(A)T\|U\|_{C([0,T]\times\R)}.
\end{equation}
From the first equation of problem \eqref{uODE} with $(v,g,h)=(v_i^\ast,g_i^\ast,h_i^\ast)$ we obtain
\begin{equation}\label{Uc1c2}
U_t(t,x^\ast)+(d+a)U(t,x^\ast)
=dc_1(t,x^\ast)+c_2(t,x^\ast),
\ 0<t\leq T,
\end{equation}
where
\begin{align*}
&c_1(t,x^\ast)
=\displaystyle\int_{g_1^\ast(t)}^{h_1^\ast(t)}J(x^\ast-y)u^\ast_1(t,y)dy
-\displaystyle\int_{g_2^\ast(t)}^{h_2^\ast(t)}J(x^\ast-y)u^\ast_2(t,y)dy,\\
&c_2(t,x^\ast)=cV(t,x^\ast).
\end{align*}
Since $U(0,x^\ast)=0$, we obtain from \eqref{Uc1c2} that
\[U(t^\ast,x^\ast)=e^{-(d+a)t^\ast}\displaystyle\int_0^{t^\ast}
e^{(d+a)\tau}[dc_1(\tau,x^\ast)+c_2(\tau,x^\ast)]d\tau.\]
For $0<t\leq T$ we have
\begin{align*}
|c_1(t,x^*)|=\ &\left|\displaystyle\int_{g_1^\ast(t)}^{h_1^\ast(t)}
J(x^*-y)u^\ast_1(t,y)dy
-\displaystyle\int_{g_2^\ast(t)}^{h_2^\ast(t)}
J(x^*-y)u^\ast_2(t,y)dy\right|\\
=\ &\left|\displaystyle\int_{g_1^\ast(t)}^{h_1^\ast(t)}J(x^*-y)
U(t,y)dy
+\left(\displaystyle\int_{g_1^\ast(t)}^{h_1^\ast(t)}
-\displaystyle\int_{g_2^\ast(t)}^{h_2^\ast(t)}\right)J(x^*-y)
u^\ast_2(t,y)dy\right|\\
\leq\ &\|U\|_{C([0,T]\times \mathbb R)}
+\left(\displaystyle\int_{g_m(t)}^{g_M(t)}
+\displaystyle\int_{h_m(t)}^{h_M(t)}\right)J(x^*-y)
|u^\ast_2(t,y)|dy\\
\leq\ &\|U\|_{([0,T]\times \mathbb R))}
+A\|J\|_\infty\left[\|g_1^\ast-g_2^\ast\|_{C([0,T])}
+\|h_1^\ast-h_2^\ast\|_{C([0,T])}\right].
\end{align*}
Therefore
\begin{equation}\label{UI}
\begin{split}
|U(t^\ast,x^\ast)|\leq\ &
C_1T\left[\|U\|_{C([0,T]\times \mathbb R))}
+\|V\|_{C([0,T]\times \mathbb R)}\right.\\
&\ \ \ \ \ \ \ \ \ +\|g_1^\ast-g_2^\ast\|_{C([0,T])}
+\|h_1^\ast-h_2^\ast\|_{C([0,T])}\Big],
\end{split}
\end{equation}
where $C_1$ depends only on $(d,c,A,J)$.
\smallskip

\underline{Case 2}: $x^\ast\in(h_0,h_m(t^\ast))$.

In this case, there exist $t_1^\ast,t_2^\ast\in(0,t^\ast)$ such that $h_1^\ast(t_1^\ast)=h_2^\ast(t_2^\ast)=x^\ast$. Without loss of generality, we may assume that $0\leq t_1^\ast\leq t_2^\ast$.

We first prove that
\begin{equation}\label{t12*}
t_2^*-t_1^*\leq \frac 1{\mu\sigma_0}\|h_1^*-h_2^*\|_{C([0,T])}.
\end{equation}
A proof is needed only if $t_1^*<t_2^*$. In such a case, from
 \[
 \frac{h_1^\ast(t_2^*)-h_1^\ast(t_1^\ast)}{t_2^*-t_1^\ast}\geq\mu\sigma_0
 \]
  we obtain
\begin{align*}
0<t_2^*-t_1^*
\leq\ &(\mu\sigma_0)^{-1}
[h_1^\ast(t_2^*)-h_1^\ast(t_1^\ast)]\\
=\ &(\mu\sigma_0)^{-1}
[h_1^\ast(t_2^*)-h_2^\ast(t_2^*)]\\
\leq\ &(\mu\sigma_0)^{-1}\|h_1^\ast-h_2^\ast\|_{C([0,T])}.
\end{align*}

We next estimate $\widetilde V(t^*,x^*)=(\widetilde{v}^\ast_1-\widetilde{v}^\ast_2)(t^\ast,x^\ast)$. Clearly
\begin{align*}
|\widetilde V(t^\ast,x^\ast)|
=\ &e^{-bt^\ast}
\left|\widetilde{v}^\ast_1(t_2^\ast,x^\ast)e^{bt_2^\ast}
+\displaystyle\int_{t_2^\ast}^{t^\ast}
e^{b\tau}\left(G(u^\ast_1)-G(u^\ast_2)\right)(\tau,x^\ast)d\tau\right|\\
\leq\ &|\widetilde{v}^\ast_1(t_2^\ast,x^\ast)|
+T\rho(A)\|U\|_{C([0,T]\times\R)}.
\end{align*}
By the first equation of \eqref{uODE} and $u_1^\ast(t_1^\ast,x^\ast)=0$, we have, for $t\in[t_1^\ast,t_2^\ast]$,
\begin{equation}\label{u1C2}
u_1^\ast(t,x^\ast)=
\displaystyle\int_{t_1^\ast}^{t}\left[
d\displaystyle\int_{g_1^\ast(\tau)}^{h_1^\ast(\tau)}J(x^\ast-y)
u_1^\ast(\tau,y)dy-du_1^\ast-au_1^\ast+cv_1^\ast\right]d\tau
\leq C_2\left(t-t_1^\ast\right),
\end{equation}
with $C_2$ depending only on $(A,B,a,c,d)$.
We have, using \eqref{t12*} and \eqref{u1C2},
\begin{align*}
0\leq \widetilde{v}^\ast_1(t_2^\ast,x^\ast)
=\ &e^{-bt_2^\ast}\displaystyle\int_{t_1^\ast}^{t_2^\ast}
e^{b\tau}G(u^\ast_1)(\tau,x^\ast)d\tau
\leq\ e^{-bt_2^\ast}\displaystyle\int_{t_1^\ast}^{t_2^\ast}
e^{bt_2^\ast}G'(0)u^\ast_1(\tau,x^\ast)d\tau\\
\leq & \ G'(0)C_2T(t_2^*-t_1^*)
\leq  \frac{G'(0)C_2T}{\mu\sigma_0}\|h_1^*-h_2^*\|_{C([0,T])}.
\end{align*}
Therefore
\begin{equation}\label{VII}
|\widetilde V(t^\ast,x^\ast)|
\leq \frac{C_2G'(0)T}{\mu\sigma_0}
\|h_1^*-h_2^*\|_{C([0,T])}+T\rho(A)\|U\|_{C([0,T]\times\R)}.
\end{equation}

Since \eqref{Uc1c2} still holds true for $t\in [t^*_2, T]$, we have
\begin{equation}\label{U-solution}
U(t^\ast,x^\ast)=e^{-(d+a)(t^\ast-t_2^\ast)}\left\{U(t_2^\ast,x^\ast)+
\displaystyle\int_{t_2^\ast}^{t^\ast}e^{(d+a)(\tau-t_2^\ast)}
[dc_1(\tau,x^\ast)+c_2(\tau,x^\ast)]d\tau\right\}.
\end{equation}
Due to $u_2^\ast(t_2^\ast,x^\ast)=0$, we have
$U(t_2^\ast,x^\ast)=u_1^\ast(t_2^\ast,x^\ast)$. Making use of \eqref{t12*} and \eqref{u1C2} we thus obtain
\[
|U(t_2^\ast,x^\ast)|\leq C_2\left(t_2^\ast-t_1^\ast\right)\leq \tilde C_2\|h_1^\ast-h_2^\ast\|_{C([0,T])},
\]
with $\tilde C_2$ depending on $C_2$ and $(\mu,\sigma_0)$. Substituting this into \eqref{U-solution}, and recalling our earlier estimates on $c_1$ and $c_2$, we obtain
\begin{equation}\label{UII}
\begin{split}
|U(t^\ast,x^\ast)|\leq\ &\tilde C_2\|h_1^\ast-h_2^\ast\|_{C([0,T])}+
TC_1\Big[\|U\|_{C([0,T]\times \mathbb R)}
+\|V\|_{C([0,T]\times \mathbb R)}\\
&\hspace{4cm}+\|g_1^\ast-g_2^\ast\|_{C([0,T])}
 +\|h_1^\ast-h_2^\ast\|_{C([0,T])}\Big].
\end{split}
\end{equation}
\smallskip

\underline{Case 3}: $x^\ast\in[h_m(t^\ast),h_M(t^\ast)]$.

Without loss of generality, we assume that $h_2^\ast(t^\ast)<h_1^\ast(t^\ast)$. In this case, there exists $t_1^\ast\in(0,t^\ast)$ such that $h_1^\ast(t_1^\ast)=x^\ast$. Then
\begin{align*}
&\widetilde{v}_2^\ast(t^\ast,x^\ast)=
u_2^\ast(t^\ast,x^\ast)=0,\ u_1^\ast(t_1^\ast,x^\ast)=0,\\
&h_m(t^\ast)=h_2^\ast(t^\ast),\
h_M(t^\ast)=h_1^\ast(t^\ast),\ h_1^\ast(t_1^\ast)=x^\ast\geq h_m(t^\ast)=h_2^\ast(t^\ast).
\end{align*}
Hence,
\[\widetilde{V}(t^\ast,x^\ast)=\widetilde{v}_1^\ast(t^\ast,x^\ast),\ U(t^\ast,x^\ast)=u_1^\ast(t^\ast,x^\ast).\]

We first prove
\begin{equation}\label{tt1-ast}
t^\ast-t_1^\ast\leq\frac{1}{\mu\sigma_0}\|h_1^\ast-h_2^\ast\|_{C([0,T])}.
\end{equation}
In fact,
\begin{align*}
0<t^\ast-t_1^\ast\leq\ &\frac{1}{\mu\sigma_0}[h_1^\ast(t^\ast)-h_1^\ast(t_1^\ast)]\leq \frac{1}{\mu\sigma_0}[h_1^\ast(t^\ast)-h_2^\ast(t^\ast)]
\leq\frac{1}{\mu\sigma_0}\|h_1^\ast-h_2^\ast\|_{C([0,T])}.
\end{align*}

By the first equation of \eqref{uODE} and $u_1^\ast(t_1^\ast,x^\ast)=0$, we have, for $t\in[t_1^\ast,t^\ast]$,
\begin{equation}\label{u1C3}
u_1^\ast(t,x^\ast)=
\displaystyle\int_{t_1^\ast}^{t}\left[
d\displaystyle\int_{g_1^\ast(\tau)}^{h_1^\ast(\tau)}J(x^\ast-y)
u_1^\ast(\tau,y)dy-du_1^\ast-au_1^\ast+cv_1^\ast\right]d\tau
\leq C_3\left(t-t_1^\ast\right),
\end{equation}
with $C_3$ depending only on $(A,B,a,c,d)$.
We next estimate
\begin{align*}
\widetilde{v}^\ast_1(t^\ast,x^\ast)
=\ &e^{-bt^\ast}\displaystyle\int_{t_1^\ast}^{t^\ast}
e^{bs}G(u^\ast_1)(s,x^\ast)ds
\leq G'(0)C_3T(t^\ast-t_1^\ast)
\leq \frac{G'(0)C_3T}{\mu\sigma_0}\|h_1^*-h_2^*\|_{C([0,T])}.
\end{align*}
Therefore,
\begin{equation}\label{VIII}
|\widetilde{V}(t^\ast,x^\ast)|\leq
\tilde C_3T\|h_1^\ast-h_2^\ast\|_{C([0,T])},
\end{equation}
with $\tilde C_3$ depending on $C_3$ and $(\mu,\sigma_0)$.

By \eqref{u1C3} and \eqref{tt1-ast}, we have
\begin{align*}
u_1^\ast(t^\ast,x^\ast)
\leq\ &C_3(t^\ast-t_1^\ast)
\leq(\mu\sigma_0)^{-1}C_3
\|h_1^\ast-h_2^\ast\|_{C([0,T])}.
\end{align*}
Thus, we obtain
\begin{equation}\label{UIII}
|U(t^\ast,x^\ast)|\leq (\mu\sigma_0)^{-1}C_3\|h_1^\ast-h_2^\ast\|_{C([0,T])}.
\end{equation}

Without loss of generality we may assume $T\leq1$.
Then the inequalities \eqref{UI}, \eqref{UII} and \eqref{UIII} yield, for the case $x^{\ast}\in[-h_{0},h_M(t^{\ast}))$,
 \begin{equation*}\label{U-C4}
\begin{split}
|U(t^\ast,x^\ast)|\leq\ &
C_4\left[T\|U\|_{C([0,T]\times\R)}
+\|V\|_{C([0,T]\times\R)}+\|g_1^\ast-g_2^\ast\|_{C([0,T])}
+\|h_1^\ast-h_2^\ast\|_{C([0,T])}\right].
\end{split}
\end{equation*}
Here $C_4$ does not depend on $T$ and $(t^*, x^*)$. When $x^{\ast}\in[g_{m}(t^{\ast}),-h_{0})$, we can similarly show that
this inequality  still holds.  Thus, since $|U(t^\ast,x^\ast)|=0$ for $0\leq t^\ast\leq T$ and $x^\ast\in\R\backslash[g_m(t^\ast),h_M(t^\ast)]$, it follows that
\begin{equation}\label{U-C6}
\begin{split}
\|U\|_{C([0,T]\times\R)}\leq\ &
C_4\left[T\|U\|_{C([0,T]\times\R)}
+\|V\|_{C([0,T]\times\R)}\right.\\
&\ \ \ \ \ \left.+\|g_1^\ast-g_2^\ast\|_{C([0,T])}
+\|h_1^\ast-h_2^\ast\|_{C([0,T])}\right].
\end{split}
\end{equation}
Hence, if $C_4T<\frac{1}{2}$, then we have
\begin{equation}\label{U-2C6}
\begin{split}
\|U\|_{C([0,T]\times\R)}\leq\ &
2C_4\left[\|V\|_{C([0,T]\times\R)}
+\|g_1^\ast-g_2^\ast\|_{C([0,T])}
+\|h_1^\ast-h_2^\ast\|_{C([0,T])}\right].
\end{split}
\end{equation}

It follows from the inequalities \eqref{VI}, \eqref{VII}, \eqref{VIII} and \eqref{U-2C6} that
\begin{equation*}
\begin{split}
|\widetilde{V}(t^\ast,x^\ast)|\leq\ &
C_5T\left[\|V\|_{C([0,T]\times\R)}
+\|g_1^\ast-g_2^\ast\|_{C([0,T])}
+\|h_1^\ast-h_2^\ast\|_{C([0,T])}\right]
\end{split}
\end{equation*}
for $0\leq t^\ast\leq T$ and $x^\ast\in [-h_0, h_M(t^*)]$. We can similarly show that this inequality also holds when $x^*\in [g_m(t^*), -h_0]$, and therefore \eqref{tilde-V-*} holds. This proves Claim 1.
\smallskip

From \eqref{tilde-V-*} we immediately obtain, for $0<T\leq \min\{\frac 1{2C_4}, T_0\}$,
\begin{equation}\label{V-C7}
\begin{split}
\|\widetilde{V}\|_{C([0,T]\times\R)}\leq\ &
C_5T\left[\|V\|_{C([0,T]\times\R)}
+\|g_1^\ast-g_2^\ast\|_{C([0,T])}
+\|h_1^\ast-h_2^\ast\|_{C([0,T])}\right]
\end{split}
\end{equation}

To derive \eqref{v-wast-dt}, we still need to estimate $\|\widetilde{g}^\ast_1-\widetilde{g}^\ast_2\|_{C([0,T])}
+\|\widetilde{h}^\ast_1-\widetilde{h}^\ast_2\|_{C([0,T])}$.

\smallskip

{\bf Claim 2:} For $0<T\leq \min\{\frac 1{2C_4}, T_0\}$, we have
\begin{equation}\label{g-h-*}
\begin{aligned}
&\|\widetilde{g}^\ast_1-\widetilde{g}^\ast_2\|_{C([0,T])}
+\|\widetilde{h}^\ast_1-\widetilde{h}^\ast_2\|_{C([0,T])}\\
\leq\ &C_6T\left(\|V\|_{C([0,T]\times\R)}
+\|g^\ast_1-g^\ast_2\|_{C([0,T])}+\|h^\ast_1-h^\ast_2\|_{C([0,T])}\right),
\end{aligned}
\end{equation}
where $C_6$ depends only on $(\mu,h_0,A,C_4)$.

For $0\leq t\leq T$,
\begin{align*}
&\widetilde{g}_i^\ast(t)=-h_0-\mu\displaystyle\int_0^t
\displaystyle\int_{g_i^\ast(\tau)}^{h_i^\ast(\tau)}
\displaystyle\int_{-\infty}^{g_i^\ast(\tau)}
J(x-y)u_i^\ast(\tau,x)dydxd\tau,\\
&\widetilde{h}_i^\ast(t)=h_0+\mu\displaystyle\int_0^t
\displaystyle\int_{g_i^\ast(\tau)}^{h_i^\ast(\tau)}
\displaystyle\int_{h_i^\ast(\tau)}^{+\infty}
J(x-y)u_i^\ast(\tau,x)dydxd\tau,\ i=1,2.
\end{align*}
Moreover, it follows from the definition of $\Sigma_T$ that
\[h_i^\ast(t)-g_i^\ast(t)\leq 2h_0+\frac{\epsilon_0}{4}\leq3h_0.\]
By direct calculations, we obtain
\begin{align*}
&|\widetilde{h}^\ast_1(t)-\widetilde{h}^\ast_2(t)|\\
\leq\ &\mu\displaystyle\int_0^t\left|
\displaystyle\int_{g_1^\ast(\tau)}^{h_1^\ast(\tau)}
\displaystyle\int_{h_1^\ast(\tau)}^{+\infty}
J(x-y)u_1^\ast(\tau,x)dydx
-\displaystyle\int_{g_2^\ast(\tau)}^{h_2^\ast(\tau)}
\displaystyle\int_{h_2^\ast(\tau)}^{+\infty}
J(x-y)u_2^\ast(\tau,x)dydx\right|d\tau\\
\leq\ &\mu\displaystyle\int_0^t
\displaystyle\int_{g_1^\ast(\tau)}^{h_1^\ast(\tau)}
\displaystyle\int_{h_1^\ast(\tau)}^{+\infty}
J(x-y)|U(\tau,x)|dydxd\tau\\
&+\mu\displaystyle\int_0^t\left|
\left(\displaystyle\int_{g_1^\ast(\tau)}^{h_1^\ast(\tau)}
\displaystyle\int_{h_1^\ast(\tau)}^{+\infty}
-\displaystyle\int_{g_2^\ast(\tau)}^{h_2^\ast(\tau)}
\displaystyle\int_{h_2^\ast(\tau)}^{+\infty}\right)
J(x-y)u_2^\ast(\tau,x)dydx\right|d\tau\\
\leq\ &3h_0\mu T
\|U\|_{C([0,T]\times\R)}\\
&+\mu\displaystyle\int_0^t\left|\left(
\displaystyle\int_{g_{m}(\tau)}^{g_{M}(\tau)}
\displaystyle\int_{h_{m}(\tau)}^{+\infty}
+\displaystyle\int_{h_{m}(\tau)}^{h_{M}(\tau)}
\displaystyle\int_{h_{m}(\tau)}^{+\infty}
+\displaystyle\int_{g_{m}(\tau)}^{h_{M}(\tau)}
\displaystyle\int_{h_{m}(\tau)}^{h_{M}(\tau)}\right)
J(x-y)u_2^\ast(\tau,x)dydx\right|d\tau\\
\leq\ &3h_0\mu T
\|U\|_{C([0,T]\times\R)}+\mu T
\left(A\|g^\ast_1-g^\ast_2\|_{C([0,T])}
+2A\|h^\ast_1-h^\ast_2\|_{C([0,T])}\right).
\end{align*}
Similarly,
\begin{align*}
&|\widetilde{g}^\ast_1(t)-\widetilde{g}^\ast_2(t)|
\leq3h_0\mu T
\|U\|_{C([0,T]\times\R)}+\mu T
\left(2A\|g^\ast_1-g^\ast_2\|_{C([0,T])}
+A\|h^\ast_1-h^\ast_2\|_{C([0,T])}\right).
\end{align*}
Combining these estimates with \eqref{U-2C6}, we immediately obtain \eqref{g-h-*}. This proves Claim 2.

\smallskip

 From \eqref{V-C7} and \eqref{g-h-*}, we deduce
\begin{align*}
&\|\widetilde{V}\|_{C([0,T]\times\R)}
+\|\widetilde{g}^\ast_1-\widetilde{g}^\ast_2\|_{C([0,T])}
+\|\widetilde{h}^\ast_1-\widetilde{h}^\ast_2\|_{C([0,T])}\\
\leq\ &(C_5+C_6)T\left(\|V\|_{C([0,T]\times\R)}
+\|g^\ast_1-g^\ast_2\|_{C([0,T])}
+\|h^\ast_1-h^\ast_2\|_{C([0,T])}\right).
\end{align*}
Therefore, if we choose $\widetilde{T}$ such that
\[0<\widetilde{T}\leq\min
\left\{T_0,\ 1,\ \frac{1}{2C_4},\ \frac{1}{2(C_5+C_6)}\right\},\]
then for any $T\in (0, \widetilde T]$, \eqref{v-wast-dt} holds with $\delta=1/2$, and so
$\mathcal{F}$ is a contraction mapping on $\Sigma_T$. Hence $\mathcal{F}$ has a unique fixed point $(v,g,h)$ in $\Sigma_T$, which gives a  nonnegative solution $(u,v,g,h)$ of \eqref{FB} for $t\in (0, T]$.
\smallskip

\textbf{Step 3:} We show that the solution $(u,v,g,h)$ of \eqref{FB} for $t\in (0, T]$ is the unique nonnegative solution of \eqref{FB} for $t\in (0, T]$.

Let $(\overline{u},\overline{v},\overline{g},\overline{h})$ be an arbitrary solution of problem \eqref{FB} for $t\in (0, T]$.
Since $(v,g,h)$ is the unique fixed point of $\mathcal F$ in $\Sigma_T$, the uniqueness conclusion will follow if we can show
$(\overline{v},\overline{g},\overline{h})\in\Sigma_T$. We first show
\begin{equation}\label{u<=A}
\overline{u}(t,x)\leq A \text{\ for\ } (t,x)\in \overline{D_{\widetilde{T}}^{g,h}}.
\end{equation}
It suffices to show that the above inequality holds with $A$ replaced by $A+\epsilon$ for any given $\epsilon>0$. Suppose this is not true. Due to $\overline{u}(0,x)<A+\epsilon=:A_{\epsilon}$, there exist $t^\ast\in(0,\widetilde{T}]$ and $x^\ast\in(\overline g(t^\ast),\overline h(t^\ast))$ such that
\[\overline{u}(t^\ast,x^\ast)=A_{\epsilon},\ \overline{u}_t(t^\ast,x^\ast)\geq0,\]
and
\[\overline{u}(t,x)<A_{\epsilon}
\text{\ for\ } t\in[0,t^\ast),\ x\in [\overline g(t),\overline h(t)].\]
Define
\begin{equation*}
\overline t_*:=
\begin{cases}
t_{x^*}^{\overline g},& \mbox{ if } x^*\in[\overline g(t^*),-h_{0}) \text{\ and\ } x^*=\overline g(t_{x^*}^{\overline g}),\\
0,& \mbox{ if }  x^*\in[-h_{0},h_{0}],\\
t_{x^*}^{\overline h},& \mbox{ if }  x^*\in(h_{0},\overline h(t^*)] \text{\ and\ }
x^*=\overline h(t_{x^*}^{\overline h}).
\end{cases}
\end{equation*}
Then $\overline{v}(t, x^*)$ solves the ODE problem
\begin{equation}\label{v}
\begin{cases}
v'=-bv+{G}(\overline{u})\leq -bv-G(A_\epsilon),&\overline t_*<t\leq t^\ast,\\
v(\overline t_*)=\hat v_0(x^*).
\end{cases}
\end{equation}
 By a simple comparison argument we obtain
\[
\overline{v}(t,x^*)\leq B_{\epsilon}
:=\max\left\{\frac{G(A_\epsilon)}{b},\|v_0\|_\infty\right\}
\text{\ for\ } t\in[\overline t_*,t^\ast].\]
Since $\frac{G(z)}{z}$ is decreasing, and $A_\epsilon>A\geq K_1$, $\frac{G(K_1)}{K_1}=\frac{ab}{c}$, it is easy to check that $B_{\epsilon}\leq\frac{a}{c}A_{\epsilon}$.
It follows that $-a\overline{u}(t^\ast,x^\ast)+c\overline{v}(t^\ast,x^\ast)
\leq-aA_{\epsilon}+cB_{\epsilon}\leq0$. Hence
\[0\leq\overline{u}_t(t^\ast,x^\ast)\leq d\displaystyle\int_{\overline g(t^\ast)}^{\overline h(t^\ast)}J(x^\ast-y)
\overline{u}(t^\ast,y)dy-d\overline{u}(t^\ast,x^\ast).\]
Since $\overline{u}(t^\ast,\overline g(t^\ast))=\overline{u}(t^\ast,\overline h(t^\ast))=0$, for $y\in(\overline g(t^\ast),\overline h(t^\ast))$ but close to the boundary of this interval, $\overline{u}(t^\ast,y)<A_\epsilon$. It follows that
\[dA_{\epsilon}=d\overline{u}(t^\ast,x^\ast)\leq d\displaystyle\int_{\overline g(t^\ast)}^{\overline h(t^\ast)}J(x^\ast-y)\overline{u}(t^\ast,y)dy
<dA_{\epsilon}\displaystyle\int_{\overline g(t^\ast)}^{\overline h(t^\ast)}
J(x^\ast-y)dy\leq dA_{\epsilon}.\]
This contradiction proves \eqref{u<=A}.

We may now use  \eqref{v} again but with $(t^*,x^*)$ replaced by an arbitrary $(t,x)\in {D_{{T}}^{\overline g,\overline h}}$, and $A_\epsilon$ replaced by $A$,  to deduce
\[\overline{v}(t,x)\leq B \text{\ for\ } (t,x)\in {D_{{T}}^{\overline g,\overline h}}.\]
 Hence, $\overline{v}\in\mathbb{X}^{v_0}_{{T}}$. Therefore $(\overline{v},\overline{g},\overline{h})\in\Sigma_{{T}}$ since the properties for $\overline g$ and $\overline h$ can be proved by
 the same argument as in step 3 of the proof of \cite[Theorem 2.1]{CDLL2018}.
 We have thus proved that problem \eqref{FB} has a unique nonnegative solution $(u,v,g,h)$ for $t\in (0, T]$.
 \smallskip

\textbf{Step 4:} Extension of the solution of \eqref{FB} to $t\in (0, \infty)$.

From Step 3 we see that problem \eqref{FB} has a unique solution $(u,v,g,h)$ over some time interval $(0, \widetilde{T}]$. For any $s\in(0,\widetilde{T})$, $(u(s,\cdot), v(s,\cdot))$ satisfies \eqref{initial} with $(-h_0, h_0)$ replaced by $(g(s), h(s))$.
This implies that we can treat $u(s,x)$ and $v(s,x)$ as the initial functions and use the above arguments to extend the solution from $t=s$ to $t=T'\geq\widetilde{T}$. Suppose that $(0,\widehat{T})$ is the maximal interval that the solution $(u,v,g,h)$ of \eqref{FB} can be defined through this extension process. We will show that $\widehat{T}=\infty$. Otherwise $\widehat{T}\in(0,\infty)$ and we will derive a contradiction.

Firstly, we can similarly show, as above, that
\[
0\leq u\leq A,\; 0\leq v\leq B \mbox{ for } (t,x)\in {D_{\widehat{T}}^{ g,h}}.
\]
For $t\in(0,\widehat{T})$, since
\begin{equation}\label{gh}
\begin{cases}
&g'(t)=-\mu\displaystyle\int_{g(t)}^{h(t)}\displaystyle\int_{-\infty}^{g(t)}
J(x-y)u(t,x)dydx,\\
&h'(t)=\mu\displaystyle\int_{g(t)}^{h(t)}\displaystyle\int_{h(t)}^{+\infty}
J(x-y)u(t,x)dydx,
\end{cases}
\end{equation}
we have
\[[h(t)-g(t)]'\leq\mu A[h(t)-g(t)].\]
Then
\[h(t)-g(t)\leq2h_0e^{\mu At} \text{\ for\ } t\in(0,\widehat{T}).\]
Since $g(t)$ and $h(t)$ are monotone functions in $(0,\widehat{T})$, we can define
\[g(\widehat{T})=\lim\limits_{t\rightarrow\widehat{T}}g(t),\
h(\widehat{T})=\lim\limits_{t\rightarrow\widehat{T}}h(t)
\text{\ with\ }
h(\widehat{T})-g(\widehat{T})\leq2h_0e^{\mu A\widehat{T}}.\]
Denote
\[\Omega_{\widehat{T}}
=\{(t,x):t\in(0,\widehat{T}],x\in(g(t),h(t))\}.\]
By \eqref{gh} and $0\leq u(t,x)\leq A$ in $\Omega_{\widehat{T}}$, we have $g',h'\in L^\infty([0,\widehat{T}))$. Together with $g(T)$ and $h(T)$ defined above, we see that $g,h\in C([0,\widehat{T}])$. It is easy to see that the right-hand sides of the first and second equations in \eqref{FB} belong to $L^\infty(\Omega_{\widehat{T}})$.
It follows that $u_t,v_t\in L^\infty(\Omega_{\widehat{T}})$. From this fact we easily see that\[u(\widehat{T},x)=\lim\limits_{t\rightarrow\widehat{T}}u(t,x)
\text{\ and\ }
v(\widehat{T},x)=\lim\limits_{t\rightarrow\widehat{T}}v(t,x)\]
exist for each $x\in(g(\widehat{T}),h(\widehat{T}))$.

We show next that $u(\widehat{T},\cdot),v(\widehat{T},\cdot)$ are continuous for $x\in(g(\widehat{T}),h(\widehat{T}))$.
For any given $\varepsilon_0>0$ small, we let
\[\theta_1=\sup\left\{\theta:
[g(\widehat{T})+\varepsilon_0,h(\widehat{T})-\varepsilon_0]
\subset(g(\widehat{T}-\theta),h(\widehat{T}-\theta))\right\}.\]
For any $\varepsilon_1>0$, we take $\theta_2=\min\left\{\frac{\varepsilon_1}{4C},\theta_1\right\}$. For any $x\in[g(\widehat{T})+\varepsilon_0,h(\widehat{T})-\varepsilon_0]$, we have
\[|u(\widehat{T},x)-u(\widehat{T}-\theta_2,x)|
=\left|\displaystyle\int_{\widehat{T}-\theta_2}^{\widehat{T}}
u_\tau(\tau,x)d\tau\right|
\leq C\theta_2\leq\frac{\varepsilon_1}{4}.\]
By the continuity of $u(\widehat{T}-\theta_2,x)$ in $x\in(g(\widehat{T}-\theta_2),h(\widehat{T}-\theta_2))$, we have $u(\widehat{T}-\theta_2,x)$ is continuous uniformly in $[g(\widehat{T})+\varepsilon_0,h(\widehat{T})-\varepsilon_0]$. Hence, for the above $\varepsilon_1>0$  there exists $\delta_1$ such that, when $|x-y|<\delta_1$ and $x,y\in[g(\widehat{T})+\varepsilon_0,h(\widehat{T})-\varepsilon_0]$, we have \[|u(\widehat{T}-\theta_2,x)-u(\widehat{T}-\theta_2,y)|
<\frac{\varepsilon_1}{2}.\]
Hence, for such $x$ and $y$,
\begin{align*}
&|u(\widehat{T},x)-u(\widehat{T},y)|\\
=\ &|u(\widehat{T},x)-u(\widehat{T}-\theta_2,x)
+u(\widehat{T}-\theta_2,x)-u(\widehat{T}-\theta_2,y)
+u(\widehat{T}-\theta_2,y)-u(\widehat{T},y)|\\
\leq\ &\varepsilon_1.
\end{align*}
This proves that $u(\widehat{T},\cdot)$ is continuous in $(g(\widehat{T}),h(\widehat{T}))$. Similarly, $v(\widehat{T},\cdot)$ is also continuous in $(g(\widehat{T}),h(\widehat{T}))$.

To show
$u\in C(\overline{\Omega}_{\widehat{T}})$ and
$v\in C(\overline{\Omega}_{\widehat{T}})$, it remains to prove that
\[\mbox{$u(t,x),\, v(t,x)\rightarrow0$ as $(t,x)\rightarrow(\widehat{T},g(\widehat{T}))$ and $(t,x)\rightarrow(\widehat{T},h(\widehat{T}))$. }
\]
We only prove the former case as the other case can be shown similarly. Noting that $t_x\nearrow\widehat{T}$ as $x\searrow g(\widehat{T})$, we have
\begin{align*}
|u(t,x)|\leq\ &\left|\displaystyle\int_{t_x}^{t}\left[
d\displaystyle\int_{g(\tau)}^{h(\tau)}J(x-y)u(\tau,y)dy-du
-au+cv\right]d\tau\right|\\
\leq\ &(t-t_x)[(2d+a)A+cB]
\rightarrow0,
\end{align*}
and
\[
0\leq v(t,x)=e^{-dt}\int_{t_x}^te^{b\tau} G(u)d\tau\leq (t-t_x)G(A)\to 0
\]
as $(t,x)\rightarrow(\widehat{T},g(\widehat{T}))$. We have shown that
$u\in C(\overline{\Omega}_{\widehat{T}})$ and
$v\in C(\overline{\Omega}_{\widehat{T}})$ and $(u,v,g,h)$ satisfies \eqref{FB} for $t\in[0,\widehat{T}]$. As before, we can use a simple comparison argument to show that  $u$ and $v$ are positive in $\Omega_{\widehat T}$. Thus we can regard $(u(\widehat{T},x),v(\widehat{T},x))$ as the initial function and extend the solution of \eqref{FB} to some $(0,T)$ with $T>\widehat{T}$. This contradicts the definition of $\widehat{T}$. Therefore we must have $\widehat{T}=\infty$.
$\hfill \Box$

\section{Long-time dynamical behaviour}

In this section, we prove Theorems 1.2 and 1.3 by a series of lemmas. Throughout this section, we always assume that $J$ satisfies {\bf (J)}, and $G$ satisfies {\bf (G1)} and {\bf (G2)}.
We start with two comparison results.

\begin{lemma}\label{MP-fixed}
Let $h_0,T>0$ and $\Omega_0:=[0,T]\times[-h_0,h_0]$. Suppose that $\alpha\in L^\infty(\Omega_0)$ is nonnegative, and $(u(t,x),v(t,x))$ as well as $(u_t(t,x),v_t(t,x))$ are continuous in $\Omega_0$ and satisfy
\begin{equation*}
\begin{cases}
u_{t}\geq d\displaystyle\int_{-h_0}^{h_0}J(x-y)u(t,y)dy-du
-au+cv,&0<t\leq T,\ x\in[-h_0,h_0],\\
v_{t}\geq -bv+\alpha(t,x)u,&0<t\leq T,\ x\in[-h_0,h_0],\\
u(0,x)\geq0,\ v(0,x)\geq0,&x\in[-h_0,h_0].
\end{cases}
\end{equation*}
Then $u(t,x),v(t,x)\geq 0$ for all $0\leq t\leq T$ and
$-h_0\leq x\leq h_0$.
\end{lemma}
\begin{proof}
This follows from a simple variation of the argument in the proof of Lemma 2.1. We omit the details.
\end{proof}

\begin{lemma}\label{comparison}
For $T\in(0,+\infty)$, suppose that $\overline{g},\overline{h}\in C([0,T])$, $\overline{u},\overline{v}\in C(\overline{D_T^{\overline{g},\overline{h}}})$, $\overline{u},\overline{v}\geq0$. If $(\overline{u},\overline{v},\overline{g},\overline{h})$ satisfies
\begin{equation*}
\begin{cases}
\overline{u}_{t}\geq d\displaystyle\int_{\overline{g}(t)}^{\overline{h}(t)}
J(x-y)\overline{u}(t,y)dy-d\overline{u}-a\overline{u}+c\overline{v},
&0<t\leq T,\ x\in(\overline{g}(t),\overline{h}(t)),\\
\overline{v}_{t}\geq -b\overline{v}+G(\overline{u}),
&0<t\leq T,\ x\in(\overline{g}(t),\overline{h}(t)),\\
\overline{g}'(t)\leq-\mu
\displaystyle\int_{\overline{g}(t)}^{\overline{h}(t)}
\displaystyle\int_{-\infty}^{\overline{g}(t)}
J(x-y)\overline{u}(t,x)dydx,&0<t\leq T,\\
\overline{h}'(t)\geq\mu
\displaystyle\int_{\overline{g}(t)}^{\overline{h}(t)}
\displaystyle\int_{\overline{h}(t)}^{+\infty}
J(x-y)\overline{u}(t,x)dydx,&0<t\leq T,\\
\overline{u}(0,x)\geq u_{0}(x),\
\overline{v}(0,x)\geq v_{0}(x),&x\in[-h_{0},h_{0}],\\
\overline{g}(0)\leq-h_0,\ \overline{h}(0)\geq h_{0},&
\end{cases}
\end{equation*}
then the unique solution $(u,v,g,h)$ of \eqref{FB} satisfies
\[u(t,x)\leq\overline{u}(t,x),\ v(t,x)\leq\overline{v}(t,x),\
g(t)\geq\overline{g}(t),\ h(t)\leq\overline{h}(t)
\text{\ for\ } 0<t\leq T, g(t)\leq x\leq h(t).\]
\end{lemma}

\begin{proof}
By {\bf (G1)}, we have $G(\overline{u})=G'(\xi)\overline{u}$ with $\xi=\xi(t,x)\in\left(0,\overline{u}(t,x)\right]$. First of all, thanks to \eqref{initial} and Lemma \ref{MP},
we have $\overline{u}>0$ for $0<t\leq T,\ \overline{g}(t)<x<\overline{h}(t)$, and thus both $\overline{h}$ and $-\overline{g}$ are strictly increasing.

For small $\epsilon>0$, let $(u_\epsilon,v_\epsilon, g_\epsilon,h_\epsilon)$ denote the unique solution of \eqref{FB} with $h_0$ replaced
by $h^\epsilon_0:= h_0(1-\epsilon)$, $\mu$ replaced
by $\mu_\epsilon:= \mu(1-\epsilon)$, and $(u_0, v_0)$ replaced by $(u^\epsilon_0, v^\epsilon_0)$ which satisfies
\[
\mbox{$0< u^\epsilon_0(x)<u_0(x)$,\; $0<v^\epsilon_0(x)<v_0(x)$ in $(-h^\epsilon_0,h^\epsilon_0)$, $u^\epsilon_0(\pm h_0^\epsilon)=v^\epsilon_0(\pm h_0^\epsilon)=0$},
\] and
\[\mbox{$\Big(u^\epsilon_0\big(\frac{h_0}{h^\epsilon_0}x\big), v^\epsilon_0\big(\frac{h_0}{h^\epsilon_0}x\big)\Big)\rightarrow (u_0(x), v_0(x))$
as $\epsilon\rightarrow0$ in the $C([-h_0,h_0])$ norm.}
\]

We claim that $h_\epsilon(t)<\overline{h}(t)$ and $g_\epsilon(t)>\overline{g}(t)$ for all $t\in(0,T]$. Clearly, these hold true for small $t>0$. Suppose that there exists $t_1\leq T$ such that
\[h_\epsilon(t)<\overline{h}(t),\ g_\epsilon(t)>\overline{g}(t) \text{\ for\ } t\in(0,t_1) \text{\ and\ } [h_\epsilon(t_1)-\overline{h}(t_1)][g_\epsilon(t_1)-\overline{g}(t_1)]=0.\]
Without loss of generality, we may assume that
\[h_\epsilon(t_1)=\overline{h}(t_1),\ g_\epsilon(t_1)\geq\overline{g}(t_1).\]
We now compare $(u_\epsilon, v_\epsilon)$ and $(\overline{u},\overline v)$ over the region
\[\Omega_{\epsilon,t_1}:=\{(t, x)\in\R^2:0<t\leq t_1, g_\epsilon(t) <x<h_\epsilon(t)\}.\]
Let $w=\overline{u}-u_\epsilon$ and $z=\overline{v}-v_\epsilon$, then $(w,z)$ satisfies
\begin{equation*}
\begin{cases}
w_{t}\geq d\displaystyle\int_{g_\epsilon(t)}^{h_\epsilon(t)}J(x-y)w(t,y)dy-dw
-aw+cz,&0<t\leq t_1,\ x\in(g_\epsilon(t),h_\epsilon(t)),\\
z_{t}\geq -bz+G'(\eta)w,&0<t\leq t_1,\ x\in(g_\epsilon(t),h_\epsilon(t)),\\
w(t,x)\geq0,\ z(t,x)\geq0,&0<t\leq t_1,\ x=g_\epsilon(t) \text{\ or\ } h_\epsilon(t),\\
w(0,x)>0,\ z(0,x) > 0,&x\in[g_\epsilon(0),h_\epsilon(0)],
\end{cases}
\end{equation*}
where $\eta=\eta(t,x)$ is between $\overline{u}(t,x)$ and $u_\epsilon(t,x)$.
By Lemma \ref{MP}, it follows that $\overline{u}>u_\epsilon,\; \overline v>v_\epsilon$ in $\Omega_{\epsilon,t_1}$.

Furthermore, according to the definition of $t_1$, we have $h'_\epsilon(t_1)\geq\overline{h}'(t_1)$. Thus
\begin{align*}
0\geq&\ \overline{h}'(t_1)-h'_\epsilon(t_1)\\
\geq&\ \mu\displaystyle\int_{\overline{g}(t_1)}^{\overline{h}(t_1)}
\displaystyle\int_{\overline{h}(t_1)}^{+\infty}
J(x-y)\overline{u}(t_1,x)dydx
-\mu_\epsilon\displaystyle\int_{g_\epsilon(t_1)}^{h_\epsilon(t_1)}
\displaystyle\int_{h_\epsilon(t_1)}^{+\infty}
J(x-y)u_\epsilon(t_1,x)dydx\\
>&\ \mu_\epsilon\displaystyle\int_{g_\epsilon(t_1)}^{h_\epsilon(t_1)}
\displaystyle\int_{h_\epsilon(t_1)}^{+\infty}
J(x-y)[\overline{u}(t_1,x)-u_\epsilon(t_1,x)]dydx>0,
\end{align*}
which is a contradiction. The claim is thus proved, i.e., we always have $h_\epsilon(t)<\overline{h}(t)$ and $g_\epsilon(t)>\overline{g}(t)$ for all $t\in(0,T]$. Moreover, we also have $\overline{u}>u_\epsilon,\; \overline v>v_\epsilon$ in $\Omega_{\epsilon,T}$.

Since the unique solution of \eqref{FB} depends continuously on the parameters in \eqref{FB}, the desired result then follows by letting $\epsilon\rightarrow0$.
\end{proof}

The following result is a direct consequence of the above comparison principle, where to stress the dependence on the parameter $\mu$, we use $(u_\mu,v_\mu,g_\mu,h_\mu)$ to denote the solution of problem \eqref{FB}.

\begin{corollary}\label{mu-comparison}
If $\mu_1\leq\mu_2$, we have
$u_{\mu_1}(t,x)\leq u_{\mu_2}(t,x),\
v_{\mu_1}(t,x)\leq v_{\mu_2}(t,x),\
g_{\mu_1}(t)\geq g_{\mu_2}(t),\
h_{\mu_1}(t)\leq h_{\mu_2}(t)$ for
$0<t\leq T, \; g_{\mu_1}(t)\leq x\leq h_{\mu_1}(t)$.
\end{corollary}

It is easily seen that $h(t)$ is monotonically increasing and $g (t)$ is monotonically decreasing. Therefore
 \[
 \lim\limits_{t\rightarrow\infty}h(t)=h_{\infty}\in (h_0,+\infty],\ \
\lim\limits_{t\rightarrow\infty}g(t)=g_{\infty}\in [-\infty, -h_0)
\]
are always well-defined.
Let
\[
\theta:=\frac{c G'(0)}{b}-a=a(R_0-1).
\]

\begin{lemma}\label{C-Vanishing}
If $\theta\leq 0$, or equivalently $R_0\leq 1$,  then $h_{\infty}-g_{\infty}<\infty$.
\end{lemma}

\begin{proof}
Direct calculations yield
\begin{align*}
&\frac{d}{dt}\int_{g(t)}^{h(t)}
\left[u(t,x)+\frac{c}{b}v(t,x)\right]dx\\
=\ &\int_{g(t)}^{h(t)}\left[u_{t}(t,x)
+\frac{c}{b}v_{t}(t,x)\right]dx
+h'(t)\left[u(t,h(t))+\frac{c}{b}v(t,h(t))\right]\\
&-g'(t)\left[u(t,g(t))+\frac{c}{b}v(t,g(t))\right]\\
=\ &\int_{g(t)}^{h(t)}\left[d\int_{g(t)}^{h(t)}J(x-y)u(t,y)dy
-d\int_{\R}J(x-y)u(t,x)dy-au+\frac{c}{b}G(u)\right]dx\\
=\ &d\int_{g(t)}^{h(t)}\int_{g(t)}^{h(t)}J(y-x)u(t,y)dxdy
-d\int_{g(t)}^{h(t)}\int_{\R}J(x-y)u(t,x)dydx\\
&+\int_{g(t)}^{h(t)}\left[-au+\frac{c}{b}G(u)\right]dx\\
=\ &d\int_{g(t)}^{h(t)}\int_{g(t)}^{h(t)}J(x-y)u(t,x)dydx
-d\int_{g(t)}^{h(t)}\int_{\R}J(x-y)u(t,x)dydx\\
&+\int_{g(t)}^{h(t)}\left[-au+\frac{c}{b}G(u)\right]dx\\
=\ &-d\int_{g(t)}^{h(t)}\int_{-\infty}^{g(t)}J(x-y)u(t,x)dydx
-d\int_{g(t)}^{h(t)}\int_{h(t)}^\infty J(x-y)u(t,x)dydx\\
&+\int_{g(t)}^{h(t)}\left[-au+\frac{c}{b}G(u)\right]dx\\
=\ &-\frac{d}{\mu}(h'(t)-g'(t))
+\int_{g(t)}^{h(t)}\left[-au(t,x)+\frac{c}{b}G(u(t,x))\right]dx.
\end{align*}
Integrating from $0$ to $t$ gives
\begin{align*}
&\int_{g(t)}^{h(t)}\left[u(t,x)+\frac{c}{b}v(t,x)\right]dx\\
=\ &\int_{g(0)}^{h(0)}
\left[u(0,x)+\frac{c}{b}v(0,x)\right]dx
+\frac{d}{\mu}(h(0)-g(0))
-\frac{d}{\mu}(h(t)-g(t))\\
&+\int_{0}^{t}\int_{g(s)}^{h(s)}
\left[-au(s,x)+\frac{c}{b}G(u(s,x))\right]dxds.
\end{align*}
Since $\frac{G(z)}{z}\leq G'(0)$ by the monotonicity of $\frac{G(z)}{z}$, it follows from $\theta\leq0$ that
\[-au(s,x)+\frac{c}{b}G(u(s,x))
\leq-au(s,x)+\frac{cG'(0)}{b}u(s,x)\leq 0.\]
Then
\[\frac{d}{\mu}(h(t)-g(t))\leq\int_{-h_0}^{h_0}
\left[u_0(x)+\frac{c}{b}v_0(x)\right]dx
+2\frac{d}{\mu}h_0,\]
which gives that $h_{\infty}-g_{\infty}<\infty$ by letting $t\rightarrow\infty$.
\end{proof}

We define the operator $\mathcal{L}_{\Omega}+\beta : C(\overline{\Omega})\rightarrow C(\overline{\Omega})$ by
\[(\mathcal{L}_{\Omega}+\beta)[\phi](x)=:
d\displaystyle\int_{\Omega}J(x-y)\phi(y)dy
-d\phi(x)+\beta(x)\phi(x),\]
where $\Omega$ is an open bounded interval in $\R$,  and $\beta\in C(\overline{\Omega})$. The generalized principal eigenvalue of $\mathcal{L}_{\Omega}+\beta$ is given by
\[\lambda_{p}(\mathcal{L}_{\Omega}+\beta)
=:\inf\left\{\lambda\in\mathbb{R}:
(\mathcal{L}_{\Omega}+\beta)[\phi]\leq\lambda\phi \text{\ in\ } \Omega \text{\ for some\ } \phi\in C(\overline{\Omega}),\ \phi>0\right\}.\]

To find the criteria for vanishing and spreading of \eqref{FB} for the case $\theta>0$, or equivalently $R_0>1$, we first consider the following fixed boundary problem, with $[l_1, l_2]$  a finite interval in $\mathbb R$,
\begin{equation}\label{wz-t-Omega}
\begin{cases}
w_t=d\displaystyle\int_{l_1}^{l_2}J(x-y)w(t,y)dy-dw-aw+cz,
&t>0,\ x\in[l_1,l_2],\\
z_t=-bz+G(w),&t>0,\ x\in[l_1,l_2],\\
w(0,x)=w_0(x),\ z(0,x)=z_0(x),&x\in[l_1,l_2].
\end{cases}
\end{equation}
The corresponding stationary problem of \eqref{wz-t-Omega} is
\begin{equation}\label{wz-Omega}
\begin{cases}
d\displaystyle\int_{l_1}^{l_2}J(x-y)w(y)dy-dw-aw+cz=0,&x\in[l_1,l_2],\\
-bz+G(w)=0,&x\in[l_1,l_2].
\end{cases}
\end{equation}

Clearly problem \eqref{wz-Omega} is equivalent to the single equation
\begin{equation}\label{w-Omega}
d\displaystyle\int_{l_1}^{l_2}J(x-y)w(y)dy-dw-aw+\frac{cG(w)}{b}=0,\
x\in[l_1,l_2].
\end{equation}
Let us note that, due to {\bf (G1)-(G2)}, when $\theta>0$,
\[
f(w):=-aw+\frac{c G (w)}{b}
\]
is a Fisher-KPP nonlinear function, namely it satisfies conditions {\bf (f3)-(f4)} in \cite{CDLL2018}. Let us also recall from \cite{CDLL2018} that,
\[
\lim_{l_2-l_1\to +\infty}\lambda_p(\mathcal L_{(l_1, l_2)}+\theta)=\theta,\; \lim_{l_2-l_1\to 0}\lambda_p(\mathcal L_{(l_1, l_2)}+\theta)=\theta-d.
\]
For \eqref{w-Omega}, by \cite{BZ2007,CDLL2018} we have the following result.
\begin{lemma}\label{U-l1l2}
If $\lambda_p(\mathcal{L}_{(l_1,l_2)}+\theta)>0$, then problem \eqref{w-Omega} admits a unique positive solution $W\in C([l_1,l_2])$, and
\[\lim\limits_{-l_1,l_2\rightarrow+\infty}W(x)=K_1
\text{\ locally uniformly in\ } \R.\]
Moreover, if $\lambda_p(\mathcal{L}_{(l_1,l_2)}+\theta)\leq0$, then any nonnegative uniformly bounded solution of \eqref{w-Omega} is identically zero.
\end{lemma}

From Lemma \ref{U-l1l2}, we obtain directly
\begin{corollary}\label{UV-l1l2}
If $\lambda_p(\mathcal{L}_{(l_1,l_2)}+\theta)>0$, then problem \eqref{wz-Omega} admits a unique positive solution $(W(x),Z(x))$ with $Z(x)=G(W(x))/b$, and
\[\lim\limits_{-l_1,l_2\rightarrow+\infty}
(W(x),Z(x))=(K_1,K_2)
\text{\ locally uniformly in\ } \R.\]
Moreover, if $\lambda_p(\mathcal{L}_{(l_1,l_2)}+\theta)\leq0$, then any nonnegative uniformly bounded solution of \eqref{wz-Omega} is identically zero.
\end{corollary}

\begin{lemma}\label{wz-t-well-posed}
  If the initial functions $w_0,\ z_0\in C([l_1,l_2])$, are nonnegative,  and $w_0,\ z_0\not\equiv0$, then \eqref{wz-t-Omega} has a unique positive solution
  $(w(t,x), z(t,x))$ defined for all $t>0$.
 \end{lemma}
 \begin{proof}
 This can be proved by arguments used in the proof of Theorem 1.1, but the proof is much more simpler, since the boundaries are fixed. We omit the details.
 \end{proof}
 The long-time dynamical behaviour of \eqref{wz-t-Omega} is determined by the following result.
\begin{lemma}\label{uv-longtimebehavior}
For $(w_0, z_0)$ as in Lemma \ref{wz-t-well-posed}, the following conclusions hold for the unique solution $(w,z)$ of \eqref{wz-t-Omega}:
\begin{itemize}
\item[{\rm(i)}] If $\lambda_p(\mathcal{L}_{(l_1,l_2)}+\theta)>0$, then
\[\lim\limits_{t\rightarrow+\infty}(w(t,x),z(t,x))=(W(x),Z(x))
\text{\ uniformly in\ } [l_1,l_2].\]
\item[{\rm(ii)}] If $\lambda_p(\mathcal{L}_{(l_1,l_2)}+\theta)\leq0$, then
\[\lim\limits_{t\rightarrow+\infty}(w(t,x),z(t,x))=(0,0)
\text{\ uniformly in\ } [l_1,l_2].\]
\end{itemize}
\end{lemma}

\begin{proof}
(i) It follows from Proposition 3.4 of \cite{CDLL2018} that $\lambda_p\leq\theta=\frac{cG'(0)}{b}-a$. Then we have
\[\frac{cG'(0)}{b}>\frac{cG'(0)}{b}-a\geq\lambda_p>\frac{\lambda_p}{4}.\]
Thus, $\frac{G'(0)}{b}-\frac{\lambda_p}{4c}>0$.
Let $\phi$ be the positive normalized eigenfunction corresponding to $\lambda_p$, namely, $\|\phi\|_{\infty}=1$ and
\begin{equation*}
d\displaystyle\int_{l_1}^{l_2}J(x-y)\phi(y)dy
-d\phi(x)-a\phi+\frac{cG'(0)}{b}\phi=\lambda_p\phi,\ x\in[l_1,l_2].
\end{equation*}
Define, for small $\epsilon>0$,
\[\underline{w}^\ast(x)=\epsilon\phi(x),\
\underline{z}^\ast(x)
=\left(\frac{G'(0)}{b}-\frac{\lambda_p}{4c}\right)\epsilon\phi(x)
\text{\ for\ } x\in[l_1,l_2].\]
Direct computations yield
\begin{align*}
&d\displaystyle\int_{l_1}^{l_2}J(x-y)\underline{w}^\ast(y)dy
-d\underline{w}^\ast-a\underline{w}^\ast+c\underline{z}^\ast\\
=\ &\epsilon\left[d\displaystyle\int_{l_1}^{l_2}J(x-y)\phi(y)dy
-d\phi-a\phi+\frac{cG'(0)}{b}\phi
-\frac{\lambda_p}{4}\phi\right]
=\frac{3}{4}\epsilon\lambda_p\phi
\end{align*}
and
\begin{align*}
&-b\underline{z}^\ast+G(\underline{w}^\ast)
=\epsilon\phi\left(\frac{b\lambda_p}{4c}-G'(0)+G'(\eta)\right),
\end{align*}
where $\eta\in(0,\underline{w}^\ast)$. Hence we can choose $\epsilon$ small enough such that $(\underline{w}^\ast(x),\underline{z}^\ast(x))$
is a lower solution of \eqref{wz-Omega}.

Since $\frac{G(K_1)}{K_1}=\frac{ab}{c}$ and $\frac{G(z)}{z}$ is decreasing, we can choose $M_1$ and $M_2$ such that
\[M_1>\max\{K_1,\|w_0\|_\infty\},\ M_2>\|z_0\|_\infty,\ G(M_1)<bM_2<\frac{ab}{c}M_1.\]
Then it is easy to check that $(\overline{w}^\ast(x),\overline{z}^\ast(x))=(M_1,M_2)$ is an upper solution of \eqref{wz-Omega}. Moreover, we can further guarantee that
\[\mbox{$(\overline{w}^\ast(x),\overline{z}^\ast(x))\geq (w(1,x), z(1,x))\geq (\underline{w}^\ast(x),\underline{z}^\ast(x))$ componentwisely.}
\]

Let $(\underline{w}(t,x),\underline{z}(t,x))$ be the solution of
\begin{equation}\label{wzl1l2}
\begin{cases}
w_t=d\displaystyle\int_{l_1}^{l_2}J(x-y)w(t,y)dy-dw-aw+cz,
&t>0,\ x\in[l_1,l_2],\\
z_t=-bz+G(w),&t>0,\ x\in[l_1,l_2],\\
w(0,x)=\underline{w}^\ast(x),\ z(0,x)=\underline{z}^\ast(x),&x\in[l_1,l_2],
\end{cases}
\end{equation}
and let $(\overline{w}(t,x),\overline{z}(t,x))$ be the solution of \eqref{wz-t-Omega} with $(w_0(x),z_0(x))$ replaced by $(\overline{w}^\ast(x),\overline{z}^\ast(x))$. By Lemma \ref{MP-fixed} and a simple comparison argument, we have
\begin{equation}\label{w-w-w}
(\underline{w}(t,x),\underline{z}(t,x))
\leq(w(t+1,x),z(t+1,x))\leq(\overline{w}(t,x),\overline{z}(t,x))
\text{\ for\ } t\geq0 \text{\ and\ } x\in[l_1,l_2].
\end{equation}
It follows from a simple comparison argument  that $(\underline{w}(t,x),\underline{z}(t,x))$ is a non-decreasing function pair in $t$, and $(\overline{w}(t,x),\overline{z}(t,x))$ is a non-increasing function pair in $t$. Therefore,
\[
\mbox{$(0,0)<(\underline{w}(t,x),\underline{z}(t,x))\leq(\overline{w}(t,x),\overline{z}(t,x))\leq(M_1,M_2)$ for $t\geq0$ and $x\in[l_1,l_2]$.}
\]
 Hence, there exists $(\underline{W}(x),\underline{Z}(x))$ such that
\begin{equation*}
\lim\limits_{t\rightarrow+\infty}
(\underline{w}(t,x),\underline{z}(t,x))= (\underline{W}(x),\underline{Z}(x))
\text{\ for all\ } x\in[l_1,l_2].
\end{equation*}

Next, we show that $(\underline{W}(x),\underline{Z}(x))$ is a solution of \eqref{wz-Omega}. It follows from the first equation of \eqref{wzl1l2} that,
for any $t,s >0$ and $x\in (l_1, l_2)$,
\begin{align*}
\underline{w}(t+s,x)-\underline{w}(t,x)=\
&\displaystyle\int_0^s\left[
d\displaystyle\int_{l_1}^{l_2}J(x-y)\underline{w}(t+\tau,y)dy
-d\underline{w}(t+\tau,x)\right.\\
&\hspace{3.5cm} -a\underline{w}(t+\tau,x)+c\underline{z}(t+\tau,x)\Big]d\tau.
\end{align*}
Similarly,
\begin{align*}
\underline{z}(t+s,x)-\underline{z}(t,x)=\ &\displaystyle\int_0^s
\left[-b\underline{z}(t+\tau,x)+G(\underline{w}(t+\tau,x))\right]d\tau.
\end{align*}
Letting $t\rightarrow\infty$, by the Lebesgue dominated convergence theorem, we obtain
\[
0=s\Big[d\displaystyle\int_{l_1}^{l_2}J(x-y)\underline{W}(y)dy
-d\underline{W}(x)-a\underline{W}(x)+c\underline{Z}(x)\Big],
\]
\[ 0=s
\left[-b\underline{Z}(x)+G(\underline{W}(x))\right].
\]
Since $s>0$ is arbitrary, it follows immediately that $(\underline{W}(x),\underline{Z}(x))$ is a positive solution of \eqref{wz-Omega}. Hence this is a continuous function pair, and we may use Dini's theorem to conclude that
\begin{equation}\label{w-W-under}
\lim\limits_{t\rightarrow+\infty}
(\underline{w}(t,x),\underline{z}(t,x))= (\underline{W}(x),\underline{Z}(x))
\text{\ uniformly in\ } [l_1,l_2].
\end{equation}

Similarly, there exists $(\overline{W}(x),\overline{Z}(x))$ such that
\begin{equation}\label{w-W-over}
\lim\limits_{t\rightarrow+\infty}
(\overline{w}(t,x),\overline{z}(t,x))= (\overline{W}(x),\overline{Z}(x))
\text{\ uniformly in\ } [l_1,l_2],
\end{equation}
and $(\overline{W}(x),\overline{Z}(x))$ is a positive solution of \eqref{wz-Omega}. Since \eqref{wz-Omega} has a unique positive solution $(W(x), Z(x))$, we must have
\begin{equation*}
(\underline{W}(x),\underline{Z}(x))=(\overline{W}(x),\overline{Z}(x))
=(W(x),Z(x)).
\end{equation*}
It follows from this fact, \eqref{w-w-w}, \eqref{w-W-under} and \eqref{w-W-over} that
\[\lim\limits_{t\rightarrow+\infty}(w(t,x),z(t,x))=(W(x),Z(x))
\text{\ uniformly in\ } [l_1,l_2].\]

(ii) In this case we can  construct $(\overline{w}(t,x),\overline{z}(t,x))$ as above, and obtain
\begin{equation*}
(0,0)\leq(w(t,x),z(t,x))\leq(\overline{w}(t,x),\overline{z}(t,x))
\text{\ for\ } t\geq0 \text{\ and\ } x\in[l_1,l_2],
\end{equation*}
\begin{equation*}
\lim\limits_{t\rightarrow+\infty}
(\overline{w}(t,x),\overline{z}(t,x))= (W(x),Z(x))
\text{\ for \ } x\in [l_1,l_2],
\end{equation*}
and $(W(x),Z(x))$ is a nonnegative solution of \eqref{wz-Omega}. Since $\lambda_p(\mathcal{L}_{(l_1,l_2)}+\theta)\leq0$, it follows from Corollary \ref{UV-l1l2} that  $(W(x),Z(x))\equiv(0,0)$.
Hence, we have
\[\lim\limits_{t\rightarrow+\infty}(w(t,x),z(t,x))=(0,0)
\text{\ uniformly in\ } [l_1,l_2].\]
The proof is complete.
\end{proof}

\begin{lemma}\label{uv-0-lam<=0}
 Let $(u,v,g,h)$ be the unique solution of \eqref{FB}. If $h_\infty-g_\infty<\infty$, then
\begin{equation}\label{uv-0}
\lim_{t\to\infty}\,\max_{g(t)\leq x\leq h(t)}u(t,x)=
\lim_{t\to\infty}\,\max_{g(t)\leq x\leq h(t)}v(t,x)=0,
\end{equation}
and
\begin{equation}\label{lam-p-L}
\lambda_p(\mathcal{L}_{(g_{\infty},h_{\infty})}+\theta)\leq0.
\end{equation}
\end{lemma}

\begin{proof}
We first prove \eqref{lam-p-L}.
Suppose on the contrary that $\lambda_p(\mathcal{L}_{(g_{\infty},h_{\infty})}+\theta)>0$. Then there exists small $\epsilon_1>0$ such that $\lambda_p(\mathcal{L}_{(g_{\infty}+\epsilon,h_{\infty}-\epsilon)}
+\theta)>0$ for $\epsilon\in(0,\epsilon_1)$. Moreover, for such $\epsilon$, there exists $T_\epsilon$ such that
\[g(t)<g_\infty+\epsilon,\ h(t)>h_\infty-\epsilon \text{\ for\ } t>T_\epsilon.\]
Consider
\begin{equation}\label{wz1}
\begin{cases}
w_t=d\displaystyle\int_{g_{\infty}+\epsilon}
^{h_{\infty}-\epsilon}J(x-y)w(t,y)dy-dw-aw+cz,
&t>T_\epsilon,\ x\in[g_{\infty}+\epsilon,h_{\infty}-\epsilon],\\
z_t=-bz+G(w),&t>T_\epsilon,\ x\in[g_{\infty}+\epsilon,h_{\infty}-\epsilon],\\
w(T_\epsilon,x)=u(T_\epsilon,x),\ z(T_\epsilon,x)=v(T_\epsilon,x),
&x\in[g_{\infty}+\epsilon,h_{\infty}-\epsilon].
\end{cases}
\end{equation}
Since $\lambda_p(\mathcal{L}_{(g_{\infty}+\epsilon,h_{\infty}-\epsilon)}
+\theta)>0$, it follows from Lemma \ref{uv-longtimebehavior} that the solution $(w(t,x),z(t,x))$ of \eqref{wz1} converges to the unique positive steady state $(W(x),Z(x))$ uniformly in $[g_{\infty}+\epsilon,h_{\infty}-\epsilon]$ as $t\rightarrow\infty$.

By Lemma \ref{MP-fixed} and a simple comparison argument, we have
\[(u(t,x),v(t,x))\geq(w(t,x),z(t,x)) \text{\ for\ } t>T_\epsilon \text{\ and\ } x\in[g_{\infty}+\epsilon,h_{\infty}-\epsilon].\]
Thus, there exists $\widehat{T}_\epsilon>T_\epsilon$ such that
\[(u(t,x),v(t,x))\geq\frac{1}{2}(W(x),Z(x)) \text{\ for\ } t>\widehat{T}_\epsilon \text{\ and\ } x\in[g_{\infty}+\epsilon,h_{\infty}-\epsilon].\]
By $J(0)>0$, there exists $\epsilon_0$ and $\delta_0$ such that $J(x)>\delta_0$ for $|x|<\epsilon_0$.
For $0<\epsilon<\min\{\epsilon_1,\epsilon_0/2\}$ and $t>\widehat{T}_\epsilon$,
\begin{align*}
h'(t)=&\ \mu\displaystyle\int_{g(t)}^{h(t)}
\displaystyle\int_{h(t)}^{+\infty}J(x-y)u(t,x)dydx\\
\geq&\ \mu\displaystyle\int_{g_{\infty}+\epsilon}
^{h_{\infty}-\epsilon}
\displaystyle\int_{h_\infty}^{+\infty}J(x-y)u(t,x)dydx\\
\geq&\ \mu\displaystyle\int_{h_{\infty}-\frac{\epsilon_0}{2}}
^{h_{\infty}-\epsilon}
\displaystyle\int_{h_\infty}^{h_\infty+\frac{\epsilon_0}{2}}
\delta_0\frac{1}{2}W(x)dydx:=\xi_0>0.
\end{align*}
 This implies that $h_\infty-g_\infty=+\infty$, which is a contradiction to our assumption $h_\infty-g_\infty<\infty$. Hence, we must have
\[\lambda_p(\mathcal{L}_{(g_{\infty},h_{\infty})}+\theta)\leq0.\]

Now we are ready to show \eqref{uv-0}. Let $(\overline{u},\overline{v})$ denote the unique solution of
\begin{equation}\label{uv-over}
\begin{cases}
\overline{u}_t=d\displaystyle\int_{g_{\infty}}
^{h_{\infty}}J(x-y)\overline{u}(t,y)dy-d\overline{u}
-a\overline{u}+c\overline{v},
&t>0,\ x\in[g_{\infty},h_{\infty}],\\
\overline{v}_t=-b\overline{v}+G(\overline{u}),&t>0,\ x\in[g_{\infty},h_{\infty}],\\
\overline{u}(0,x)=\widehat{u}_0(x),\ \overline{v}(0,x)=\widehat{v}_0(x),
&x\in[g_{\infty},h_{\infty}],
\end{cases}
\end{equation}
where
\begin{equation*}
\widehat{u}_{0}(x):=
\begin{cases}
u_{0}(x),&x\in[-h_{0},h_{0}],\\
0,&x\not\in[-h_{0},h_{0}]
\end{cases}
\text{\ and\ }
\widehat{v}_{0}(x):=
\begin{cases}
v_{0}(x),&x\in[-h_{0},h_{0}],\\
0,&x\not\in[-h_{0},h_{0}].
\end{cases}
\end{equation*}
By Lemma \ref{comparison}, we have $u(t,x)\leq\overline{u}(t,x)$ and $v(t,x)\leq\overline{v}(t,x)$ for $t>0$ and $x\in[g(t),h(t)]$.
Since $\lambda_p(\mathcal{L}_{(g_{\infty},h_{\infty})}+\theta)\leq0$,
Lemma \ref{uv-longtimebehavior} implies that \[\lim\limits_{t\rightarrow\infty}(\overline{u}(t,x),
\overline{v}(t,x))=(0,0) \text{\ uniformly in\ } x\in[g_\infty,h_\infty].\]
Hence \eqref{uv-0} holds.
\end{proof}

\begin{lemma}\label{C-Spreading}
If $\theta\geq d$, then $h_{\infty}-g_{\infty}=\infty$.
\end{lemma}

\begin{proof}
Arguing indirectly we assume that $h_{\infty}-g_{\infty}<\infty$. Thanks to \cite[Proposition 3.4]{CDLL2018},
\[\lambda_p(\mathcal{L}_{(g_{\infty},h_{\infty})}+\theta)>0.\]
This is a contradiction to \eqref{lam-p-L}.
\end{proof}

We next consider the case
\[0<\theta<d.\]
In this case, it follows from \cite[Proposition 3.4]{CDLL2018} that there exists $l^{\ast}$ such that
\begin{equation}\label{eigen}
\begin{cases}
\lambda_{p}(\mathcal{L}_{(l_{1},l_{2})}+\theta)=0 &\mbox{ if } l_{2}-l_{1}=l^{\ast}\\
 (\l_2-l_1-l^*)\lambda_{p}(\mathcal{L}_{(l_{1},l_{2})}+\theta)>0 & \mbox{ if } l_{2}-l_{1}\in (0,+\infty)\setminus\{l^*\}
 \end{cases}
\end{equation}

\begin{lemma}\label{C-Spreading-Vanishing}
Assume that $0<\theta<d$. Then the following hold:
\begin{itemize}
\item[{\rm(i)}] If $h_{\infty}-g_{\infty}<\infty$, then $h_{\infty}-g_{\infty}\leq l^\ast$.
\item[{\rm(ii)}] If $2h_0\geq l^\ast$, then $h_{\infty}-g_{\infty}=\infty$.
\item[{\rm(iii)}] If $2h_{0}<l^\ast$, then there exist $0<\mu_0\leq \mu^0<+\infty$ such that $h_{\infty}-g_{\infty}=\infty$ for $\mu>\mu^0$ and $h_{\infty}-g_{\infty}<\infty$ for $0<\mu<\mu_0$.
\end{itemize}
\end{lemma}

\begin{proof}
(i) Arguing indirectly we assume that $h_{\infty}-g_{\infty}>l^\ast$. Since $0<\theta<d$, we have $\lambda_p(\mathcal{L}_{(g_{\infty},h_{\infty})}+\theta)>0$. This is a contradiction to \eqref{lam-p-L}.

(ii) This conclusion follows directly from (i).

(iii) By  \cite[Lemma 3.9]{DWZ2019},  there exists $\mu^0$ such that $h_{\infty}-g_{\infty}=\infty$ for $\mu>\mu^0$.
It remains to prove the conclusion for $\mu_0$.

Since $2h_{0}<l^\ast$, we have
$\lambda_{p}(\mathcal{L}_{(-h_{0},h_{0})}+\theta)<0$.
There exists some small $\varepsilon>0$ such that $h^{\ast}:= h_{0}\left(1+\varepsilon\right)$ satisfies
\[\lambda_{p}^{1}:=
\lambda_{p}(\mathcal{L}_{(-h^{\ast},h^{\ast})}+\theta)<0.\]
Let $\phi_1$ be the positive normalized eigenfunction corresponding to $\lambda_p^1$, namely, $\|\phi_1\|_{\infty}=1$ and
\begin{equation}\label{lambda-phi-i}
d\displaystyle\int_{-h^{\ast}}^{h^{\ast}}J(x-y)\phi_1(y)dy
-d\phi_1(x)-a\phi_1+\frac{cG'(0)}{b}\phi_1=\lambda_p^1\phi_1,\ x\in[-h^{\ast},h^{\ast}].
\end{equation}
Choose a positive constant $K$ large enough such that
\[K\phi_1(x)\geq u_{0}(x) \text{\ and\ }
\left(\frac{G'(0)}{b}-\frac{\lambda_p^1}{4c}\right)K\phi_1(x)\geq v_{0}(x) \text{\ for\ } x\in[-h_0,h_0].\]
Define
\begin{align*}
&\overline{h}(t)=h_{0}
\left[1+\varepsilon\left(1-e^{-\delta t}\right)\right],\ \overline{g}(t)=-\overline{h}(t),\ t\geq0,\\
&\overline{u}(t,x)=Ke^{-\delta t}\phi_{1}(x),
\ t\geq0,\ x\in[\overline{g}(t),\overline{h}(t)],\\
&\overline{v}(t,x)=\left(\frac{G'(0)}{b}-\frac{\lambda_p^1}{4c}\right)
\overline{u}(t,x),\ t\geq0,\ x\in[\overline{g}(t),\overline{h}(t)],
\end{align*}
where $\delta>0$ will be determined later. Clearly $h_{0}\leq\overline{h}(t)\leq h^{\ast}$.

Direct calculations yield
\begin{align*}
&\overline{u}_{t}-d\displaystyle\int_{\overline{g}(t)}^{\overline{h}(t)}
J(x-y)\overline{u}(t,y)dy
+d\overline{u}+a\overline{u}-c\overline{v}\\
\geq\ &Ke^{-\delta t}\left(-\delta\phi_1(x)
-d\displaystyle\int_{\overline{g}(t)}^{\overline{h}(t)}
J(x-y)\phi_1(y)dy+d\phi_1+a\phi_1
-\frac{cG'(0)}{b}\phi_1+\frac{\lambda_p^1}{4}\phi_1\right)\\
\geq\ &Ke^{-\delta t}\left(-\delta
-\frac{3\lambda_p^1}{4}\right)\phi_1(x)
\end{align*}
and
\begin{align*}
\overline{v}_{t}+b\overline{v}-G(\overline{u})
>\ &(b-\delta)\left(\frac{G'(0)}{b}-\frac{\lambda_p^1}{4c}\right)
\overline{u}-G'(0)\overline{u}\\
=\ &\left[-\delta\frac{G'(0)}{b}
-(b-\delta)\frac{\lambda_p^1}{4c}\right]\overline{u}
\end{align*}
for $t>0$ and $x\in(\overline{g}(t),\overline{h}(t))$. Since $\lambda_p^1<0$, we can choose $\delta$ small enough such that
\[-\delta-\frac{3\lambda_p^1}{4}>0 \text{\ and\ }
-\delta\frac{G'(0)}{b}-(b-\delta)\frac{\lambda_p^1}{4c}>0.\]
Moreover, $\overline{h}'(t)=h_{0}\varepsilon\delta e^{-\delta t}$ and
\begin{align*}
&\mu\displaystyle\int_{\overline{g}(t)}^{\overline{h}(t)}
\displaystyle\int_{\overline{h}(t)}^{+\infty}
J(x-y)\overline{u}(t,x)dydx\leq2\mu Ke^{-\delta t}h^\ast.
\end{align*}
If
\[\mu\leq\frac{h_{0}\varepsilon\delta}{2Kh^{\ast}}:=\mu_{0},\]
then we have
\[\overline{h}'(t)\geq\mu\displaystyle\int_{\overline{g}(t)}^{\overline{h}(t)}
\displaystyle\int_{\overline{h}(t)}^{+\infty}
J(x-y)\overline{u}(t,x)dydx.\]
Similarly, we can derive
\[\overline{g}'(t)\leq-\mu\displaystyle\int_{\overline{g}(t)}^{\overline{h}(t)}
\displaystyle\int_{-\infty}^{\overline{g}(t)}
J(x-y)\overline{u}(t,x)dydx.\]
We may now apply Lemma \ref{comparison} to obtain
\[g(t)\geq\overline{g}(t),\ h(t)\leq\overline{h}(t),\ u(t,x)\leq\overline{u}(t,x),
\ v(t,x)\leq\overline{v}(t,x) \text{\ for\ } t>0 \text{\ and\ } x\in[g(t),h(t)].\]
It follows that
$\lim\limits_{t\rightarrow\infty}(h(t)-g(t))\leq
\lim\limits_{t\rightarrow\infty}(\overline{h}(t)-\overline{g}(t))\leq 2h^{\ast}<\infty$.
\end{proof}

\begin{lemma}\label{mu}
Assume that $0<\theta<d$.
If $2h_{0}<l^\ast$, then there exist $\mu^\ast>0$ such that $h_{\infty}-g_{\infty}=\infty$ for $\mu>\mu^\ast$ and $h_{\infty}-g_{\infty}<\infty$ for $0<\mu\leq\mu^\ast$.
\end{lemma}

\begin{proof}
Define
\[\Sigma=\{\mu: \mu>0 \text{\ such that\ } h_\infty-g_\infty<+\infty\}.\]
By Lemma \ref{C-Spreading-Vanishing}, we see that $0<\sup\Sigma<+\infty$. Let $(u_\mu,v_\mu,g_\mu,h_\mu)$ denote the solution of \eqref{FB}, and set $h_{\mu,\infty}:=\lim\limits_{t\rightarrow\infty}h_\mu(t),\
g_{\mu,\infty}:=\lim\limits_{t\rightarrow\infty}g_\mu(t)$, and denote $\mu^\ast=\sup\Sigma$.

According to Corollary \ref{mu-comparison}, $u_\mu,v_\mu,-g_\mu,h_\mu$ are increasing in $\mu>0$. It follows immediately  that if $\mu_1\in\Sigma$, then $\mu\in\Sigma$ for $\mu<\mu_1$, and if $\mu_1\not\in\Sigma$, then $\mu\not\in\Sigma$ for any $\mu>\mu_1$. Hence
\begin{equation}\label{mu-Sigma}
(0,\mu^\ast)\subseteq\Sigma,\ (\mu^\ast,+\infty)\cap\Sigma=\emptyset.
\end{equation}

To complete the proof, it remains to show that $\mu^\ast\in\Sigma$. Suppose that $\mu^\ast\not\in\Sigma$. Then $h_{\mu^\ast,\infty}=-g_{\mu^\ast,\infty}=+\infty$. Thus there exists $T>0$ such that $-g_{\mu^\ast}(t)>l^\ast,\ h_{\mu^\ast}(t)>l^\ast$ for $t\geq T$. Hence there exists $\epsilon>0$ such that for $|\mu-\mu^\ast|<\epsilon$, $-g_{\mu}(T)>\frac{l^\ast}{2},\ h_{\mu}(T)>\frac{l^\ast}{2}$, which implies $\mu\not\in\Sigma$. This clearly contradicts \eqref{mu-Sigma}. Therefore $\mu^\ast\in\Sigma$.
\end{proof}

\begin{lemma}\label{iaoi} Assume that $\theta>0$. Then
$h_\infty=+\infty$ if and only if $g_\infty=-\infty$.
\end{lemma}
\begin{proof}
This follows the idea in the proof  of Lemma 3.8 in \cite{CDLL2018}. For example, if $g_\infty=-\infty$ but $h_\infty<+\infty$, then we may argue as in the proof of Lemma 3.8 here to obtain $h'(t)\geq \xi_0>0$ for all large $t$, which yields a contradiction.
\end{proof}

\begin{lemma}\label{u-v-infty} Suppose $\theta>0$.
If $h_\infty-g_\infty=+\infty$, then
\[\lim_{t\rightarrow\infty}(u(t,x),v(t,x))=(K_1,K_2)
\text{\ locally uniformly in\ } \R.\]
\end{lemma}
\begin{proof}
By Lemma \ref{iaoi}, $h_\infty-g_\infty=+\infty$ implies that $h_\infty=-g_\infty=+\infty$. Then we can choose an increasing sequence $\{t_n\}$ satisfying
\[\lim\limits_{n\rightarrow\infty}t_n=+\infty,\
\lambda_p(\mathcal{L}_{(g(t_n),h(t_n))}+\theta)>0
\text{\ for\ } n\geq1.\]
Denote $g_n=g(t_n)$ and $h_n=h(t_n)$. Let $(\underline{u}_n(t,x),\underline{v}_n(t,x))$ be the unique solution of
\begin{equation}\label{under-v}
\begin{cases}
\underline{u}_t=d\displaystyle\int_{g_n}^{h_n}J_2(x-y)\underline{u}(t,y)dy
-d\underline{u}-a\underline{u}+c\underline{v},
&t>t_n,\ x\in(g_n,h_n),\\
\underline{v}_t=-b\underline{v}+G(\underline{u}),
&t>t_n,\ x\in(g_n,h_n),\\
\underline{u}(t_n,x)=u(t_n,x),\ \underline{v}(t_n,x)=v(t_n,x),&x\in[g_n,h_n].
\end{cases}
\end{equation}
Since $\lambda_p(\mathcal{L}_{(g_n,h_n)}+\theta)>0$, it follows from Lemma \ref{uv-longtimebehavior} that \eqref{under-v} admits a unique positive steady state $(\underline{U}_n(x),\underline{V}_n(x))$ and
\[\lim\limits_{t\rightarrow\infty}
(\underline{u}_n(t,x),\underline{v}_n(t,x))
=(\underline{U}_n(x),\underline{V}_n(x))
\text{\ uniformly in\ } [g_n,h_n].\]
By Corollary \ref{UV-l1l2},
\[\lim\limits_{n\rightarrow\infty}(\underline{U}_n(x),\underline{V}_n(x))
=(K_1,K_2) \text{\ locally uniformly in\ } \R.\]
Applying Lemma \ref{MP-fixed} and a simple comparison argument, we obtain
\[(u(t,x),v(t,x))\geq(\underline{u}_n(t,x),\underline{v}_n(t,x)) \text{\ for\ } t\geq t_n \text{\ and\ } x\in[g_n,h_n].\]
Hence,
\[\liminf\limits_{t\rightarrow\infty}(u(t,x),v(t,x))\geq(K_1,K_2)
\text{\ locally uniformly in\ } \R.\]

To complete the proof, it remains to show that
\begin{equation}\label{uv<=K12}
\limsup\limits_{t\rightarrow\infty}(u(t,x),v(t,x))\leq(K_1,K_2)
\text{\ locally uniformly in\ } \R.
\end{equation}
Let $(\widehat{u}(t),\widehat{v}(t))$ be the solution of
\begin{equation}\label{uv-abcG}
\begin{cases}
u'(t)=-au(t)+cv(t),&t>0,\\
v'(t)=-bv(t)+G(u(t)),&t>0,\\
u(0)=\|u_0\|_\infty,\ v(0)=\|v_0\|_\infty.&
\end{cases}
\end{equation}
By Lemma \ref{comparison}, we have $u(t,x)\leq\widehat{u}(t)$ and $v(t,x)\leq\widehat{v}(t)$ for $t>0$ and $x\in[g(t),h(t)]$. Since $\theta>0$ and hence $R_0>1$, the unique positive equilibrium $(K_1,K_2)$  of \eqref{uv-abcG}
is globally attractive  and so $(\widehat{u}(t),\widehat{v}(t))\rightarrow(K_1,K_2)$ as $t\rightarrow\infty$, which clearly implies \eqref{uv<=K12}.
\end{proof}

Clearly Theorems \ref{dichotomy} and \ref{criteria} follow directly from Lemmas \ref{C-Vanishing}, \ref{uv-0-lam<=0}, \ref{C-Spreading}, \ref{C-Spreading-Vanishing}, \ref{mu} and \ref{u-v-infty}.

\section*{Acknowledgments}
\noindent

M. Zhao was supported by a scholarship from the China Scholarship Council. Y. Zhang was supported by NSF of China (11626072), W.-T. Li was supported by NSF of China (11731005, 11671180) and Y. Du  was supported by the Australian Research Council (DP190103757).


\end{document}